\documentclass[preprint,3p,fleqn,sort&compress,numbers]{elsarticle}
\usepackage{amsthm,amsmath,amsfonts,amssymb}
\usepackage{booktabs}
\usepackage[geometry]{ifsym}
\usepackage{multirow}
\usepackage{mathtools}
\usepackage{esint}
\usepackage{graphicx}
\usepackage{colortbl}
\usepackage{bm}
\usepackage[usenames,dvipsnames]{xcolor}
\usepackage{caption,subcaption}
\usepackage{dsfont}
\usepackage[toc,page]{appendix}
\usepackage[colorlinks]{hyperref}
\usepackage{cleveref}

\newdefinition{remark}{Remark}

\newcommand{\abs}[1]{\left|#1\right|}

\definecolor{review_color}{rgb}{0.0, 0.0, 0.0}

\crefname{appsec}{Appendix}{Appendices}

\graphicspath{{./fig/}}

\begin{document}

\definecolor{py_green}{rgb}{0.0, 0.5019607843137255, 0.0}
\definecolor{py_blue}{rgb}{0.0, 0.0, 1.0}
\definecolor{py_purple}{rgb}{0.5019607843137255, 0.0, 0.5019607843137255}
\definecolor{py_orange}{rgb}{1.0, 0.6470588235294118, 0.0}


\title{Immersed Boundary Smooth Extension: A high-order method for solving PDE on arbitrary smooth domains using Fourier spectral methods}

\author[ucd]{David~B.~Stein\corref{cor1}}
\ead{dbstein@math.ucdavis.edu}

\author[ucd]{Robert~D.~Guy}

\author[ucd]{Becca~Thomases}

\address[ucd]{Department of Mathematics, University of California, Davis, Davis, CA 95616-5270, USA}
\cortext[cor1]{Corresponding author}

\begin{abstract}
	The Immersed Boundary method is a simple, efficient, and robust numerical scheme for solving PDE in general domains, yet it only achieves first-order spatial accuracy near embedded boundaries.  In this paper, we introduce a new high-order numerical method which we call the Immersed Boundary Smooth Extension (IBSE) method.  The IBSE method achieves high-order accuracy by smoothly extending the unknown solution of the PDE from a given smooth domain to a larger computational domain, enabling the use of simple Cartesian-grid discretizations (e.g. Fourier spectral methods).  The method preserves much of the flexibility and robustness of the original IB method.  In particular, it requires minimal geometric information to describe the boundary and relies only on convolution with regularized delta-functions to communicate information between the computational grid and the boundary.  We present a fast algorithm for solving elliptic equations, which forms the basis for simple, high-order implicit-time methods for parabolic PDE and implicit-explicit methods for related nonlinear PDE.  We apply the IBSE method to solve the Poisson, heat, Burgers', and Fitzhugh-Nagumo equations, and demonstrate fourth-order pointwise convergence for Dirichlet problems and third-order pointwise convergence for Neumann problems.
\end{abstract}

\begin{keyword}
	Embedded boundary, Immersed Boundary, Fourier spectral method, Complex geometry, Partial Differential Equations, High-order
\end{keyword}

\maketitle



\section{Introduction}

The Immersed Boundary (IB) method was originally developed for the study of moving, deformable structures immersed in a fluid, and it has been widely applied to such problems since its introduction \cite{Peskin2002,Mittal2005,olson2014simulating}.  Recently, the method has been adapted to more general fluid-structure problems, including the motion of rigid bodies immersed in a fluid \cite{kallemov2015immersed} and fluid flow through a domain with either stationary boundaries or boundaries with prescribed motion \cite{Taira2007,Teran2009}.  In this broadened context, we use the term \emph{Immersed Boundary method} to refer only to methods in which ({\small\emph{i}}) the boundary is treated as a Lagrangian structure embedded in a geometrically simple domain, ({\small\emph{ii}}) the background PDE (e.g. the Navier-Stokes equations) are solved on a Cartesian grid everywhere in that domain, and ({\small\emph{iii}}) all communication between the Lagrangian structure and the underlying PDE {\color{review_color}is} mediated only by convolutions with regularized $\delta$-functions.  These methods have many desirable properties: they make use of robust and efficient Cartesian-grid methods for solving the underlying PDE, are flexible to a wide range of problems, and are simple to implement, requiring minimal geometric information and processing to describe the boundary.

The IB method belongs to the broad category of methods known as \emph{embedded boundary} (EB) methods, including the Immersed Interface \cite{li2006immersed}, Ghost Fluid \cite{Fedkiw1999}, and Volume Penalty methods \cite{Angot1999}.  These methods share a common feature: they enable solutions to PDE on nontrivial domains to be computed using efficient and robust structured-grid discretizations; yet these methods differ largely in how boundary conditions are enforced and whether or not the solution is produced in the entirety of a simple domain.  Methods which compute the solution everywhere in a $d$-dimensional rectangle admit the simplest discretizations and enable the use of high-order discretizations such as Fourier spectral methods.  Unfortunately, this simplicity comes coupled with a fundamental difficulty: the \emph{analytic} solution to these problems is rarely globally smooth on the entire domain.  Consider the one-dimensional Poisson problem $\Delta u=f$ on the periodic interval $\mathbb{T}=[0,2\pi]$ with Dirichlet boundary conditions $u(a)=u(b)=0$ for $a\neq b\in\mathbb{T}$.  Even if $f\in C^\infty(\mathbb{T})$, the solution $u$ will typically display jumps in its derivative at the values $x=a$ and $x=b$ (see \Cref{fig:extension:analytic}).  The lack of regularity in the analytic problem leads to low-order convergence in many numerical schemes, including the Immersed Boundary method.

The advantages of EB methods are substantial enough that significant effort has been expended on improving their accuracy \cite{Lai2000,Mark2008,Linnick2005,Liu2014,Xu2006,Zhong2007,Yu2007,Zhou2006,Gibou2005,Boyd2005,Bueno-Orovio2006,Lui2009,Sabetghadam2009,Albin2011,Lyon2010a,Lyon2010,Shirokoff2013}.  Two different approaches are generally taken.  The first approach involves locally altering the discretization of the PDE in the vicinity of the boundary to accommodate the lack of smoothness in the solution.  One example of this approach is the Immersed Interface method \cite{li2006immersed}.  Such approaches are particularly useful for interface problems where the solution is required on both sides of the embedded boundary.  When the solution is only needed on one side of the interface, a second approach may be taken in which variables are redefined outside of the domain of interest to obtain higher global regularity.  Improved convergence rates are achieved as a natural consequence of the properties of the discretization scheme when applied to smooth problems.  Variations of this basic idea have been used by the Fourier Continuation (FC) \cite{Lyon2010a,Lyon2010,Albin2011} and the Active Penalty (AP) \cite{Shirokoff2013} methods to provide high order solutions to PDE on general domains.

In this paper, we introduce a new method termed the \emph{Immersed Boundary Smooth Extension} (IBSE) method.  This method is a first step in remedying two deficiencies of Immersed Boundary methods:
\begin{enumerate}
	\item For generic problems, the IB method produces discretizations that are only \emph{first-order accurate} in the vicinity of domain boundaries (or fluid-structure interfaces) \cite{Beyer1992}.
	\item The IB method is only able to specify Dirichlet boundary conditions (no-slip, for fluid problems) without specialized interpolation schemes subordinate to the geometry \cite{Pacheco-Vega2007}.  Although this is not of interest when solving traditional IB fluid-interface problems, it allows the IB method to be generalized for solving other PDE (i.e. reaction-diffusion equations).
\end{enumerate}
In this paper we restrict our attention to PDE set on stationary domains, and consider only PDE without a global divergence constraint.  Direct-forcing IB methods produce solutions to PDE that are $C^0$, with jumps in the normal-derivative of the solution across the boundary, and converge at first order in the grid-spacing $\Delta x$ \cite{kallemov2015immersed,Taira2007}.  Drawing on ideas from the AP and FC methods, we use Fourier spectral methods to obtain a highly accurate discretization, while adding a volumetric forcing to non-physical portions of the computational domain to force the solution to be \emph{globally $C^k$}.  We will use the shorthand IBSE-$k$ to refer to our method when we need to explicitly denote the global regularity of the solution that is enforced.  High-order accuracy is achieved naturally, as a simple consequence of the convergence properties of the Fourier transform.

The IBSE method retains the essential robustness and simplicity of the original IB method.  All communication between the Lagrangian boundary and the underlying Cartesian grid is achieved by convolution with regularized $\delta$-functions or normal derivatives of those $\delta$-functions.  This allows an absolute minimum of geometric information to be used.  In the traditional IB method, only the position of the Lagrangian structure must be known; the IBSE method additionally requires normals to that structure and an indicator variable denoting whether Cartesian grid points lie inside or outside of the physical domain where the PDE is defined.  Additionally, since normal derivatives can be accurately approximated by convolution with derivatives of regularized $\delta$-functions, Neumann and Robin boundary conditions can be imposed in the same way that Dirichlet boundary conditions are in direct-forcing IB approaches.

In contrast to other approaches based on the smooth extension of the forcing function or the solution \cite{Lyon2010a,Lyon2010,Albin2011,Shirokoff2013,Boyd2005}, we do not extend the forcing function or solution from a previous timestep.  Instead, we \emph{smoothly extend the unknown solution} to the entire computational domain.  This approach directly enforces smoothness of the solution, and it allows for high-order \emph{implicit}-time discretizations for parabolic equations and \emph{implicit-explicit} discretizations for many nonlinear PDE.  Remarkably, the coupled problem for the solution to an elliptic PDE and that solution's smooth extension can be reduced to the solution of a relatively small dense system of equations (with size a small multiple of the number of points used to discretize the boundary), along with several FFTs.  The dense system depends only on the boundary and the discretization, and so it can be formed and prefactored to allow for efficient time-stepping of parabolic equations.

This paper is organized as follows.  In \Cref{section:poisson_methods}, we introduce the method, ignoring the details of the numerical implementation.  The fundamental contributions of this paper are contained in the modification to the IB discretization that enforces smoothness of the solution, and in the system of equations (given in \Cref{subsection:poisson_methods:system}) that allows for the simultaneous solution of the PDE along with the smooth extension of that unknown solution.  In \Cref{section:numerical_implementation}, we discuss the particular numerical implementation that we choose, detailing a fast algorithm for solving elliptic equations.  In \Cref{section:poisson_test}, we compute solutions to the Poisson problem in one and two dimensions, demonstrating high-order pointwise convergence: up to fourth-order for Dirichlet problems and third-order for Neumann problems.  In \Cref{section:heat}, we discuss discretization of the heat equation, detailing a fast algorithm for \emph{implicit} timestepping, and demonstrate fourth-order convergence in space and time.  These numerical tests include direct comparisons to the Fourier Continuation \cite{Lyon2010a} and Active Penalty \cite{Shirokoff2013} methods.  Finally, in \Cref{section:nonlinear}, we solve two nonlinear problems: a 2D Burgers' equation with homogeneous Dirichlet boundary conditions, and the Fitzhugh-Nagumo equations with homogeneous Neumann boundary conditions.  We demonstrate fourth-order convergence in space and time for Burgers' equation and third-order convergence for the Fitzhugh-Nagumo equations.



\section{Methods}
\label{section:poisson_methods}

\begin{figure}
	\centering
	\hspace*{\fill}
	\begin{subfigure}[b]{0.3\textwidth}
		\centering
		\includegraphics[width=\textwidth]{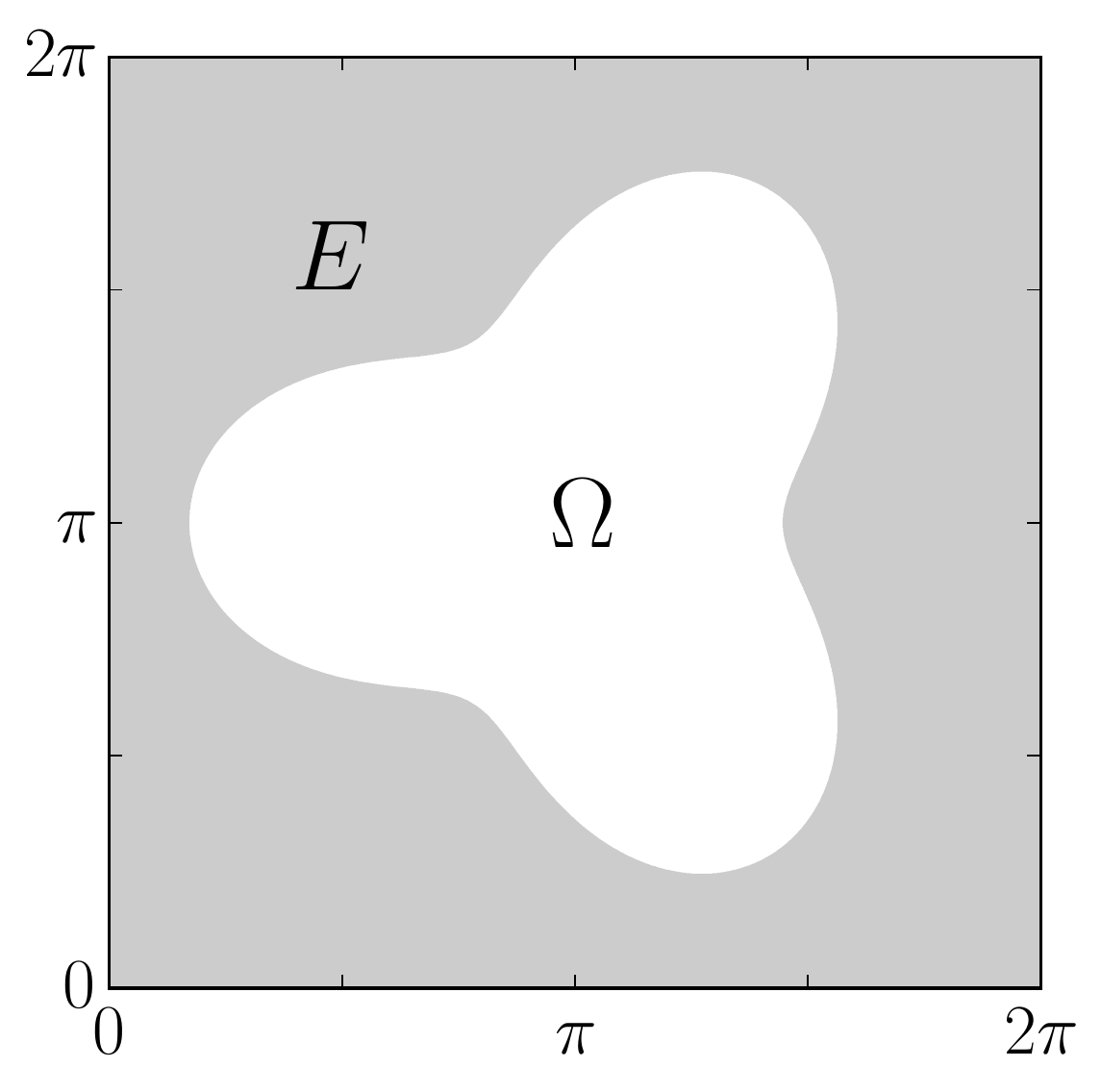}
		\subcaption{Irregular domain}
		\label{fig:domain:1}
	\end{subfigure}
	\hfill
	\begin{subfigure}[b]{0.3\textwidth}
		\centering
		\includegraphics[width=\textwidth]{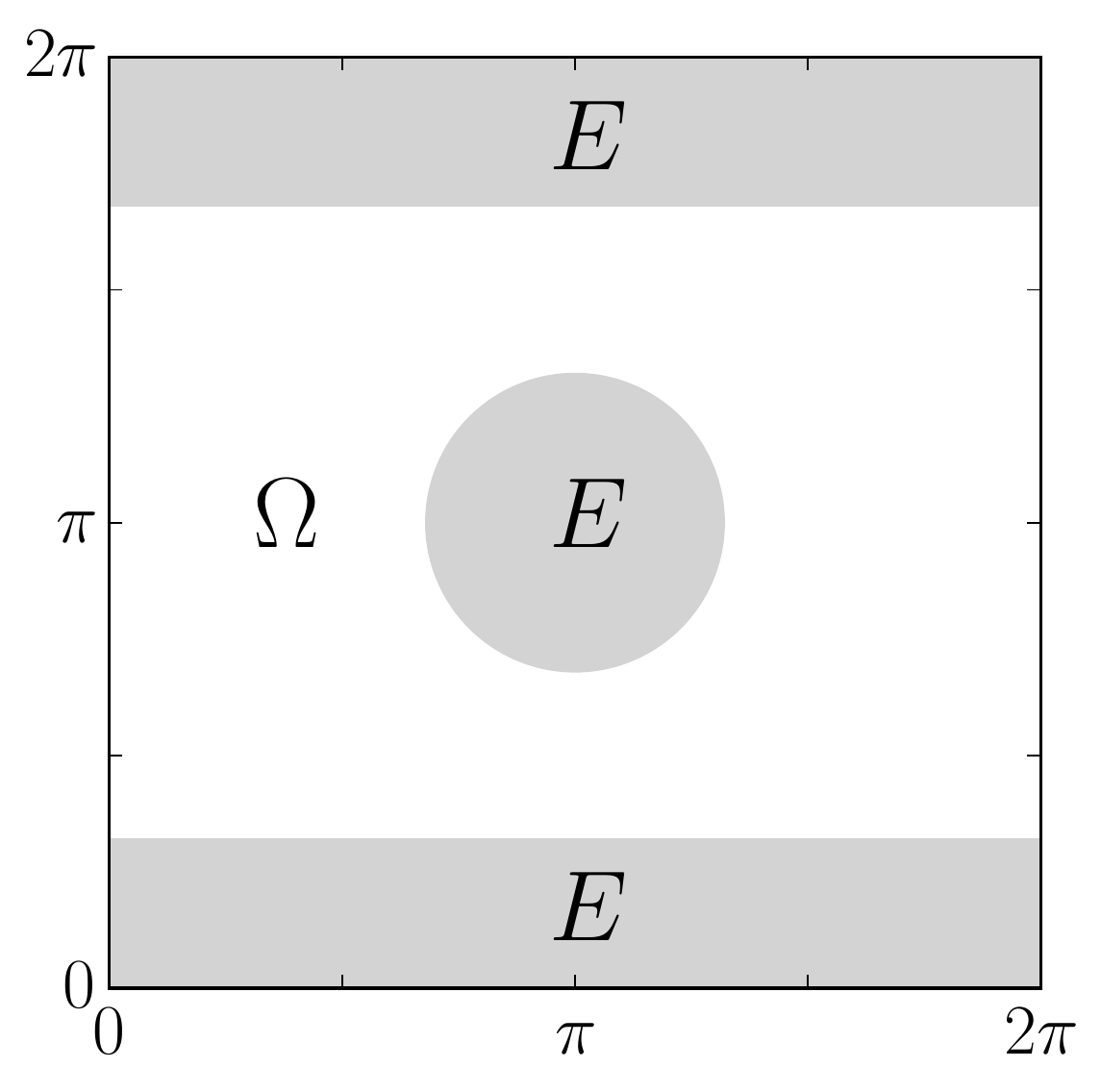}
		\subcaption{Flow around cylinder}
		\label{fig:domain:2}
	\end{subfigure}
	\hspace*{\fill}
	\caption{Two different and typical domains to solve PDE on.  The physical domain $\Omega$ is shown in white, the extension domain $E$ is shown in gray.  In both cases, we have taken the computational domain $C$ to be the 2-torus, $\mathbb{T}^2=[0,2\pi]\times[0,2\pi]$.  \Cref{fig:domain:1}, shows an irregular, connected domain $\Omega$ embedded in $C$.  \Cref{fig:domain:2} shows a domain that would be used to compute flow around a cylinder; here the domain $\Omega$ is itself periodic in one dimension, while the extension region $E$ is not connected.}
	\label{fig:domain}
\end{figure}

Let $\Omega\subset C$, $\Gamma=\partial\Omega$, and $E=C\setminus\overline{\Omega}$.  We will assume that $\Omega$ is smooth and does not self-intersect.  Two typical domains are shown in \Cref{fig:domain}.  We refer to $\Omega$ as the \emph{physical domain}, $E$ as the \emph{extension domain}, and $C$ as the \emph{computational domain}.  We first consider the Poisson problem with Dirichlet boundary conditions:
\begin{subequations}
	\label{eq:poisson}
	\begin{align}
		\Delta u	&=	f	&	&	\text{in }\Omega,	\label{eq:poisson:momentum}	\\
		u			&=	g	&	&	\text{on }\Gamma. \label{eq:poisson:boundary}
	\end{align}
\end{subequations}
For now, we assume that $f$ and $g$ are smooth ($C^\infty$) functions defined in $\Omega$ and $\Gamma$ respectively.  Let $f_e$ be a smooth extension of $f$ to $C$, that is, $f_e\in C^\infty(C)$ and $f_e(x)=f(x)$ for all $x\in\Omega$.  
The \emph{direct forcing} formulation of the Immersed Boundary method provides a way to solve \Cref{eq:poisson}:
\begin{subequations}
	\label{eq:poisson_whole_domain}
	\begin{align}
		\Delta u(x) + \int_\Gamma G(s)\delta(x-s)\,ds	&=	f_e(x)		&	&	\text{in }C,	\\
		\int_C u(x)\delta(x-s)\,dx					&=	g(s)			&	&	\text{in }\Gamma.
	\end{align}
\end{subequations}
\Cref{eq:poisson:momentum} is solved in the entire computational domain $C$, while the boundary condition is enforced by the addition of a singular force $G$ supported on the boundary that acts as a Lagrange multiplier.  This singular force term leads to jumps in the normal derivative of the solution $u$ across the boundary $\Gamma$.  This lack of smoothness restricts Immersed Boundary formulations such as \Cref{eq:poisson_whole_domain} to first-order convergence in the vicinity of the boundary unless one-sided discrete $\delta$-functions are used, which make the implementation of the method more difficult \cite{Lai2000,Beyer1992}.

\begin{figure}
	\centering
	\hspace*{\fill}
	\begin{subfigure}[b]{0.4\textwidth}
		\centering
		\includegraphics[width=\textwidth]{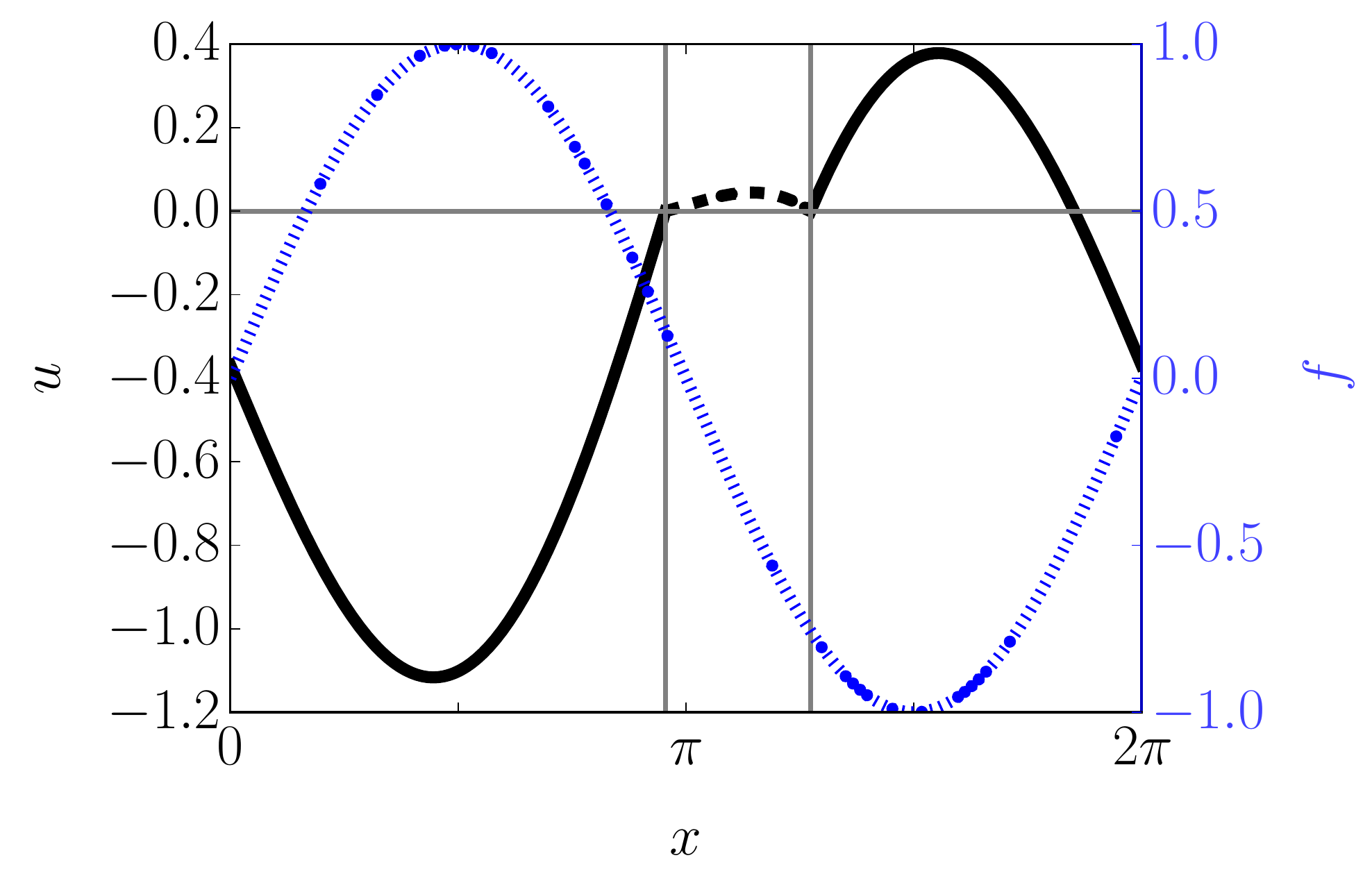}
		\subcaption{with extended forcing $f_e$}
		\label{fig:extension:analytic}
	\end{subfigure}
	\hfill
	\begin{subfigure}[b]{0.4\textwidth}
	\centering
	\includegraphics[width=\textwidth]{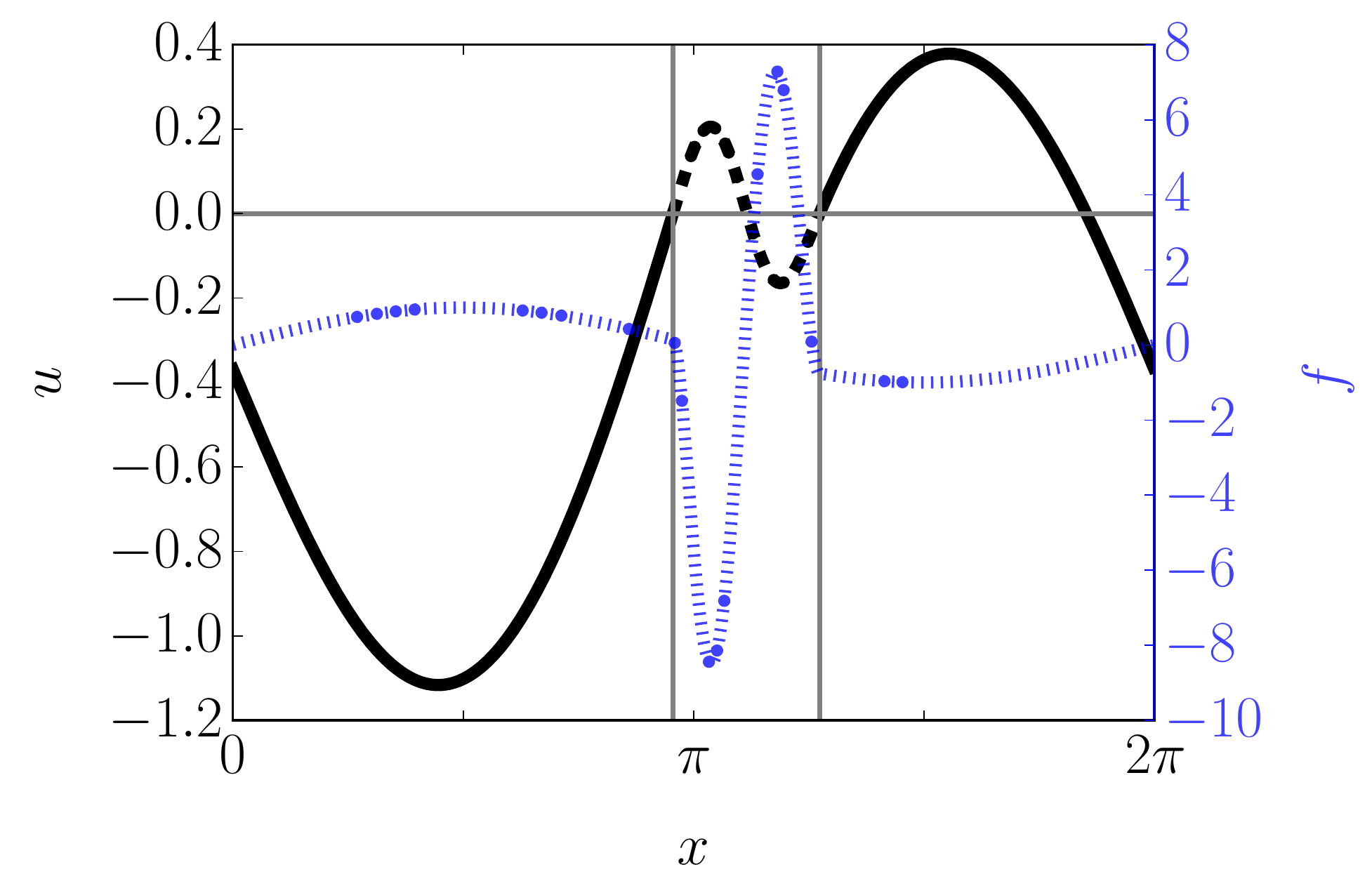}
	\subcaption{with extended forcing $\tilde f_e$}
	\label{fig:extension:smooth}
	\end{subfigure}
	\hspace*{\fill}
	\caption{The solid lines show the analytic solution $u$ to \Cref{eq:example} in the physical domain $\Omega$.  The dashed line gives the solution $u$ to the extended problem in \Cref{eq:poisson_whole_domain}, computed using the extended forcing function $f_e$ in the left figure, and $\tilde f_e$ in the right figure.  The extended forcing functions are shown as the dotted lines.  Note the different scales on the forcing functions in the two figures; the forcing is equal everywhere in the physical domain $\Omega$.}
	\label{fig:extension}
\end{figure}

We define an example one-dimensional problem to illustrate the limited regularity of solutions produced by \Cref{eq:poisson_whole_domain}.  Let the computational domain be $C=\mathbb{T}$, where the one-dimensional torus $\mathbb{T}$ is identified with the periodic interval $[0,2\pi]$, and let the physical domain be given by the interval $\Omega=\mathbb{T}\setminus[3,4]$.  The extension domain for this problem is $E=[3,4]$.  On this domain, we will solve the problem:
\begin{subequations}
	\label{eq:example}
	\begin{align}
		\Delta u	&=	\sin x	&	&\text{in }\Omega=\mathbb{T}\setminus[3,4],	\\
		u			&=	0		&	&\text{on }\Gamma=\{3,4\}.
	\end{align}
\end{subequations}
As $f$ is chosen to be $f=\sin x$ in $\Omega$, we may choose $f_e=\sin x$ for all $x\in\mathbb{T}$.  The analytic solution to this problem is shown in \Cref{fig:extension:analytic}, along with the solution to the extended problem given by \Cref{eq:poisson_whole_domain}.  Note that the solution $u$ is \emph{only continuous} despite the fact that $f_e\in C^\infty(\mathbb{T})$.  Choosing $f_e$ to be smooth leads to low regularity in the solution when the solution to the associated homogeneous problem is nontrivial, as is the case here.  Without modification, direct discretizations of this problem will exhibit slow convergence due to the limited regularity of $u$ in the vicinity of the boundary.

Rather than choose a smooth extension to the forcing function $f$, we choose an extension of the forcing function that gives a smooth solution on the entire computational domain.  Let $u_e$ be \emph{any} smooth extension of the unknown solution $u$ into $C$.  We can compute a forcing function $F_e$ associated with $u_e$:
\begin{equation}
	\label{eq:compute_forcing}
	F_e = \Delta u_e.
\end{equation}
The extended forcing function $f_e$ is then defined to be
\begin{equation}
	\label{eq:extended_forcing}
	\tilde f_e = \chi_\Omega f + \chi_E F_e,
\end{equation}
where $\chi_X$ denotes the characteristic function of the domain $X$.  \Cref{fig:extension:smooth} shows the extended forcing function $\tilde f_e$, along with the associated smooth solution $u$ to \Cref{eq:example}.  Since the solution $u$ is smooth, we expect that standard discretizations of \Cref{eq:poisson_whole_domain} will converge more rapidly when computed using the extension $\tilde f_e$ than when using the extension $f_e$.  We emphasize that the extended forcing function $\tilde f_e$ \emph{depends on the unknown solution $u$}.

This motivates us to define a reformulated version of \Cref{eq:poisson_whole_domain}, {\color{review_color}with the extended forcing function $f_e$ given by $\tilde f_e$ as defined in \Cref{eq:extended_forcing}}:
\begin{subequations}
	\label{eq:poisson_whole_domain:smooth}
	\begin{align}
		\Delta u(x) - (\chi_E\Delta\mathcal{E}_ku)(x) + \int_\Gamma G(s)\delta(x-s)\,ds	&=	\chi_\Omega f(x)	&	&	\text{in }C,	
		\label{eq:poisson_whole_domain:smooth:main} \\
		\int_C u(x)\delta(x-s)\,dx					&=	g(s)			&	&	\text{in }\Gamma.
	\end{align}
\end{subequations}
Here $\mathcal{E}_k$ is \emph{any} smooth extension operator that satisfies
\begin{subequations}
	\label{eq:extension_operator_definition}
	\begin{align}
		\mathcal{E}_k:C^0(C)\cap C^k(\Omega)	&\to	C^k(C),		\label{eq:extension_operator_definition:exten}
\\
		(\mathcal{E}_ku)(x) &= u(x)	\qquad\forall x\in\Omega.
		\label{eq:extension_operator_definition:boundary}
	\end{align}
\end{subequations}
{\color{review_color}
When restricted to the physical domain $\Omega$, \Cref{eq:poisson_whole_domain:smooth:main} reduces to $\Delta u = f$, the problem of physical interest.  When restricted to the extension domain $E$, \Cref{eq:poisson_whole_domain:smooth:main} reduces to $\Delta u = \Delta\mathcal{E}_ku$.  As long as $u\in C^0(C)$, then \Cref{eq:extension_operator_definition:boundary} ensures that $u|_\Gamma=(\mathcal{E}_ku)|_\Gamma$, and hence $u=\mathcal{E}_ku$ in $E$.  Again by \Cref{eq:extension_operator_definition:boundary}, $u=\mathcal{E}_ku$ in $\Omega$, and so $u=\mathcal{E}_ku$ in $C$.  Since $\mathcal{E}_k u\in C^k(C)$ by \Cref{eq:extension_operator_definition:exten}, then $u\in C^k(C)$.  This additional regularity of the solution $u$ allows standard discretizations of \Cref{eq:poisson_whole_domain:smooth} to converge at a faster rate than discretizations of \Cref{eq:poisson_whole_domain}.  Using $C=\mathbb{T}^d$ and a simple Fourier-spectral discretization, numerical solutions to \Cref{eq:poisson_whole_domain:smooth} should converge at $\mathcal{O}(\Delta x^{k+1})$ \cite{boyd2001chebyshev}.
}

In the remainder of this section, we lay out the remaining components needed to fully specify the IBSE method \emph{that do not depend on the particular choice of numerical implementation}.  In \Cref{subsection:poisson_methods:extension}, we discuss how to smoothly extend a known function from $\Omega$ to $C$.  In \Cref{subsection:poisson_methods:system}, we outline a system of equations that allows us to solve for $u$ and its smooth extension simultaneously.  A numerical implementation of this method will be discussed in \Cref{section:numerical_implementation}.

\subsection{Smooth extension of a known function from $\Omega$ to $C$}
\label{subsection:poisson_methods:extension}

One way to construct a smooth extension to a function is to solve a high-order PDE.  Compared to the localized and effectively one-dimensional extension strategies taken in \cite{Lyon2010a,Shirokoff2013}, the decision to extend a function by solving a fully $d$-dimensional PDE may appear to be needlessly complex.  However, we will show in \Cref{subsection:poisson_methods:system,section:numerical_implementation} that this choice leads to a robust method that is straightforward to implement and requires relatively low computational cost.

In particular, to compute a $C^k(C)$ extension to a function $v\in C^k(\Omega)$, we solve the following equation in the extension domain $E$:
\begin{subequations}
	\label{eq:extension}
	\begin{align}
		\mathcal{H}^k\xi				&=	0	&	&	\text{in }E,	\\
		\frac{\partial^j\xi}{\partial n^j}	&=	\frac{\partial^jv}{\partial n^j}	&	&\text{on }\Gamma,\ \forall\, 0\leq j\leq k.
	\end{align}
\end{subequations}
Here $\partial^j\xi/\partial n^j$ denotes the $j^\text{th}$ normal derivative of $\xi$ on the boundary $\Gamma$; a computational formula for $\partial^j\xi/\partial n^j$ is given later in \Cref{eq:normal_derivative_formula}.  A simple choice of the operator $\mathcal{H}^k$ is the polyharmonic operator $\mathcal{H}^k = \Delta^{k+1}$.  From an analytic perspective, this choice is sufficient, although we will show in \Cref{section:numerics:extension} that other choices of $\mathcal{H}^k$ will produce better results due to numerical issues.  The smooth extension $\zeta=\mathcal{E}_kv$ may then be constructed as
\begin{equation}
	\zeta = \chi_\Omega v + \chi_E\xi.
\end{equation}

\begin{remark}
	Our choice of smoothed forcing $\tilde f_e$ \emph{only depends on $\mathcal{E}_ku$ in the extension region $E$} (see \Cref{eq:extended_forcing}), so there is no need to compute the actual extension $\zeta$ since the values of $\mathcal{E}_ku$ in $\Omega$ are irrelevant.  For our purposes, we can think of an `extension' as simply a $C^k$ function in the computational domain $C$ that shares its first $k$ normal derivatives with $u$ on the boundary $\Gamma$, i.e. the function $\xi$.
\end{remark}

As with \Cref{eq:poisson}, we may solve \Cref{eq:extension} in $\Omega$ by extending the problem to all of $C$.  This can be accomplished by adding singular forces supported along the boundary:
\begin{subequations}
	\label{eq:extension_whole_domain}
	\begin{align}
		\mathcal{H}^k\xi(x) + \sum_{j=0}^{k}(-1)^j\int_\Gamma F_j(s)\frac{\partial^j\delta(x-s)}{\partial n^j}\,ds	&=	0	&	&\forall{x}\in\mathbb{T}^d,	\\
		(-1)^j\int_C\xi(x)\frac{\partial^j\delta(x-s)}{\partial n^j}\,dx	&=	\frac{\partial^j v}{\partial n^j}(s)	&	&\forall s\in\Gamma,\ 0\leq j\leq k.
		\label{eq:extension_whole_domain:boundary}
	\end{align}
\end{subequations}
The integral on the left-hand side of \Cref{eq:extension_whole_domain:boundary} is notation for the action of the distribution $\partial^j\delta/\partial n^j$ on the smooth function $\xi$.  The boundary condition given in \Cref{eq:extension_whole_domain:boundary} forces the solution $\xi$ to be $C^k$, so there is no need to alter this formulation to provide additional regularity.

For notational convenience, we define the operators $T_k$ and $T_k^*$ by:
\begin{subequations}
	\begin{align}
		T_kF(x) &= \sum_{j=0}^{k}(-1)^j\int_\Gamma F_j(s)\frac{\partial^j\delta(x-s)}{\partial n^j}\,ds,	\\
		T^*_k\xi(s) &= 
		\begin{pmatrix}
			\int_C\xi(s)\delta(x-s)\,ds							&
			\int_C\xi(s)\frac{\partial\delta(x-s)}{\partial n}	&
			\cdots									&
			\int_C\xi(s)\frac{\partial\delta^k(x-s)}{\partial n^k}
		\end{pmatrix}^\intercal.
	\end{align}
\end{subequations}
In the language of the Immersed Boundary method, $T_k^*$ is an \emph{interpolation} operator; $T_k$ is a \emph{spread} operator.  The notation is suggestive of the fact that these operators obey the adjoint property
\begin{equation}
	\left<u,T_k F\right>_C = \left<T_k^*u,F\right>_\Gamma \label{eq:adjoint_property}
\end{equation}
{\color{review_color}for all sufficiently smooth $u$ and $F$,} where the notation $\big<\cdot,\cdot\big>_X$ denotes the $L^2$ inner product on $X$.  Using the operators $T_k$ and $T^*_k$, we may represent \Cref{eq:extension_whole_domain} as
\begin{equation}
	\label{eq:extension_discretized}
	\begin{pmatrix}
		\mathcal{H}^k	&	T_k	\\
		T^*_k			&
	\end{pmatrix}
	\begin{pmatrix}
		\xi	\\	F
	\end{pmatrix}
	=
	\begin{pmatrix}
		0	\\	T_k^* v
	\end{pmatrix}.
\end{equation}

\subsection{A coupled system of equations for $u$ and its extension}
\label{subsection:poisson_methods:system}

Let the spread ($S$) and interpolation ($S^*$) operators be defined as
\begin{subequations}
	\begin{align}
		SG(x) 		&= \int_\Gamma G(s)\delta(x-s)\,ds, \\
		S^*\xi(s)	&= \int_C \xi(x)\delta(x-s)\,dx, \quad\forall s\in\Gamma.
	\end{align}
\end{subequations}
Notice that $S=T_0$ and $S^*=T_0^*$.  We can now represent \Cref{eq:poisson_whole_domain:smooth} simply as
\begin{subequations}
	\begin{align}
		\Delta u - \chi_E\Delta\xi + SG &= \chi_\Omega f,	\\
		S^*u &= g,
	\end{align}
	\label[equation]{eq:the_system:part_1}
\end{subequations}
where $\xi$ is defined by the extension equation
\begin{equation}
	\mathcal{H}^k\xi = 0,
	\label{eq:the_system:part_2}
\end{equation}
along with the constraint that
\begin{equation}
	T_k^*\xi = T_k^*u.
	\label{eq:the_system:part_3}
\end{equation}
\Cref{eq:the_system:part_1,eq:the_system:part_2,eq:the_system:part_3} can be combined into one system of equations:
\begin{equation}
	\label{eq:the_system}
	\begin{pmatrix}
		\Delta	&	-\chi_E\Delta	&		&	S	\\
				&	\mathcal{H}^k	&	T_k	&		\\
		-T_k^*	&	T_k^*			&		&		\\
		S^*		&					&		&		
	\end{pmatrix}
	\begin{pmatrix}
		u	\\	\xi	\\	F	\\	G
	\end{pmatrix}
	=
	\begin{pmatrix}
		\chi_\Omega	f	\\	0	\\	0	\\	g
	\end{pmatrix}.
\end{equation}
\Cref{eq:the_system} is the system of equations that allows the IBSE method to solve for $u$ and its extension $\xi$ simultaneously.  The remainder of this paper will be concerned with the discretization and inversion of \Cref{eq:the_system} and numerical studies demonstrating the accuracy of solutions produced by those discretizations.  We remark that \Cref{eq:the_system} is equally valid for operators of the Helmholtz type ($\Delta - \kappa\mathbb{I}$) that appear in the discretization of parabolic equations, which will be discussed in \Cref{section:heat}.



\section{Numerical implementation}
\label{section:numerical_implementation}

For concreteness, we will discuss the numerical implementation in the context of Fourier spectral methods, with the computational domain $C$ given by $C=\mathbb{T}^d=[0,2\pi]^d$.  We will discretize the domain using a regular Cartesian mesh with $n$ points $x_n$ discretizing each dimension: $\Delta x=2\pi/n$ and $x_n=n\Delta x$.  Differential operators are discretized in the usual way.

\begin{remark}
Few of the details depend upon the choice to use Fourier spectral methods, and they are chosen because of their simple implementation, computational speed, and high order of accuracy.  {\color{review_color}Inversion of the elliptic operator $\mathcal{L}$, as well as the extension operator $\mathcal{H}^k$ is also greatly simplified when using Fourier spectral methods.}  However, any discretization based on a regular Cartesian mesh could be used with minimal modifications to the method.
\end{remark}

In order to fully describe a discretization and solution strategy to \Cref{eq:the_system}, we must describe several key elements.
\begin{enumerate}
	\item In \Cref{section:numerics:operators}, we discretize the spread ($S$, $T_k$) and interpolation ($S^*$, $T_k^*$) operators.
	\item In \Cref{section:numerics:operators:delta} we define a new regularization of the $\delta$-function that is accurate and smooth.
	\item In \Cref{section:numerics:extension}, we define our extension operator $\mathcal{H}^k$ and show how this choice of $\mathcal{H}^k$ controls the numerical conditioning of the extension problem.
	\item In \Cref{section:numerics:inversion}, we describe an efficient inversion strategy for the IBSE-$k$ system given by \Cref{eq:the_system}.
	\item Finally, in \Cref{section:numerics:complexity}, we provide implementational details and briefly discuss the computational complexity of our inversion scheme.
\end{enumerate}

\subsection{Discretization of $S$, $S^*$, $T_k$, and $T_k^*$}
\label{section:numerics:operators}

Let the boundary $\Gamma$ be parametrized by the function $X(s)$.  In all examples in this manuscript, the boundary $\Gamma$ is one-dimensional and closed; the single parameter $s$ is defined on the periodic interval $[0,2\pi]$.  The spread operator $S:\Gamma\to C$ is defined as
\begin{equation}
	(S F)(x) = \int_\Gamma F(s)\delta(x-X(s))\,ds.
\end{equation}
The discrete version of this operator requires a regularized $\delta$-function and a discretization of the integral over $\Gamma$.  Construction of a regularized $\delta$-function with the properties necessary for the IBSE-$k$ method is non-trivial; we will delay discussion of this choice until \Cref{section:numerics:operators:delta} and simply denote the regularized $\delta$-function as $\tilde\delta$.  Multivariate $\delta$-functions are computed as Cartesian products of the univariate $\delta$ and also denoted by $\tilde\delta$.  Discretization of the integral over $\Gamma$ is made by choosing $n_\text{bdy}$ quadrature nodes $\tilde\Gamma=\{X_i\}_{i=1}^{n_\text{bdy}}$, equally spaced in the parameter interval $[0,2\pi]$ so that $X_i=X(s_i)$ and $s_i=(i-1)2\pi/n_\text{bdy}$.  Quadrature weights are computed at each quadrature node to be $\Delta s_i=\left\|\frac{\partial X}{\partial s}(s_i)\right\|_2$; this is a spectrally accurate quadrature rule for the integral of smooth periodic functions on $\Gamma$.  The discrete spread operator $S$ maps functions sampled at points in $\tilde\Gamma$ to $C$ by
\begin{equation}
	(SF)(x) = \sum_{i=1}^{n_\text{bdy}} F(s_i)\tilde\delta(x-X_i)\Delta s_i.
\end{equation}
We do not adopt explicit notation to distinguish between the analytic and discretized operators.  The number of points in the quadrature is chosen so that $\Delta s\approx 2\Delta x$.  This choice of node-spacing is wider than that recommended for the traditional IB method \cite{Peskin2002} but has been observed empirically to be the optimal choice in other studies of \emph{direct-forcing} IB methods \cite{kallemov2015immersed}.  The interpolation operator $S^*$ may be defined by the adjoint property $\left<u,SF\right>_C=\left<S^*u,F\right>_\Gamma$, but we note that the discrete interpolation operator $S^*$ produces a discrete function
\begin{equation}
	(S^*u)(s_k) = \int_C u(x)\delta(x-X(s_k))\,dx.
\end{equation}
Discrete integrals over $C$ are straightforward sums computed over the underlying uniform Cartesian mesh.

Normal derivatives of $\tilde\delta$ are computed by the formula \cite{john1982partial}
\begin{equation}
	\frac{\partial^j\tilde\delta}{\partial n^j} = n_{i_1}\cdots n_{i_2}\frac{\partial^j\tilde\delta}{\partial x_{i_1}\cdots \partial x_{i_j}},
	\label{eq:normal_derivative_formula}
\end{equation}
where the Einstein summation convention has been used to indicate sums over repeated indices and $\partial^j\tilde\delta/\partial x_{i_1}\cdots \partial x_{i_j}$ is computed as Cartesian products of the appropriate derivatives of the one-dimensional $\tilde\delta$.  For example, $\partial\tilde\delta/\partial n$ in two dimensions is computed as
\begin{equation}
	\frac{\partial\tilde\delta}{\partial n} = n_x \tilde\delta' \otimes \tilde\delta + n_y\tilde\delta\otimes\tilde\delta'.
\end{equation}
Analogous to the definition of $S$, we define the spread operator for the $j^\text{th}$ normal derivative, $T_{(j)}$, as
\begin{equation}
	(T_{(j)}F)(x) = (-1)^j\sum_{i=1}^{n_\text{bdy}} F(s_i)\frac{\partial^j\tilde\delta}{\partial n^j} (x-X_i)\Delta s_i
\end{equation}
and again define the interpolation operator $T^*_{(j)}$ by the adjoint property $\left<u,T_{(j)}F\right>_C=\big<T^*_{(j)}u,F\big>_\Gamma$.  Analogous to the action of $S^*$, the operator $T^*_{(j)}$ produces approximations of the $j^\text{th}$ normal derivative of a function defined on $C$ at the quadrature nodes of $\tilde\Gamma$.
Finally, the composite operator $T^*_k$ is defined by:
\begin{equation}
	T^*_k =
	\begin{pmatrix}
		T^*_{(0)}	&
		T^*_{(1)}	&
		\cdots		&
		T^*_{(k)}
	\end{pmatrix}^\intercal,
\end{equation}
while $T_k=\sum_{j=1}^kT_{(j)}$.

\subsection{Construction of a smooth discretization of $\delta$}
\label{section:numerics:operators:delta}

The IBSE-$k$ method is capable of producing solutions that converge at $\mathcal{O}(\Delta x^{k+1})$.  In order to achieve this accuracy, our choice of regularized $\delta$-function must satisfy several conditions.  For a general value of $k$:
\begin{enumerate}
	\item Its interpolation accuracy must be at least $\mathcal{O}(\Delta x^{k+1})$.
	\item It must be at least $C^k$, so that its $k^\text{th}$ normal derivative is continuous, allowing it to be used in the discretization of $T_k$ and $T^*_k$.
\end{enumerate}
In addition, we will require that the $\delta$-function has compact (and small) support so that the spread and interpolation operators $S$, $S^*$, $T_k$, and $T_k^*$ may be rapidly applied.  The commonly used `four point' $\delta$-function is $C^1$, has a support of four gridpoints, and produces approximations of $S$, $S^*$, $T_1$ and $T^*_1$ that are accurate to $\mathcal{O}(\Delta x^2)$ \cite{Griffith2005}.  The use of this $\delta$-function regularization is sufficient for the implementation of IBSE-$1$, but it is not sufficient for higher order methods, due to both its limited \emph{accuracy} and \emph{regularity}.

In this manuscript, we discretize the IBSE-$k$ method for $k=1$, $2$, and $3$, corresponding to second, third, and fourth order accuracy in $\Delta x$ for Dirichlet problems.  We therefore need a regularization of $\delta$ that is $C^3$ and has an interpolation accuracy of $\mathcal{O}(\Delta x^4)$ for smooth functions.  We are not aware of any functions with these properties currently defined in the literature, and so we construct such a function here.  For simplicity, we do not attempt to impose other conditions that are often imposed on $\delta$-functions, such as the \emph{even-odd} condition or \emph{sum of squares} condition \cite{Liu2012}.  The strategy that we follow is to choose a $\delta$-function with sufficient accuracy but limited regularity and convolve it against itself to increase its smoothness.  A similar approach was used in \cite{Yang2009} to generate a smoother $(C^2)$ version of the traditional Peskin four point $\delta$-function \cite{Peskin2002}.  In order to provide an analytic formula, a base $\delta$-function with a simple functional form must be used.  We start with the function
\begin{equation}
	\delta_{IB_4}(r) = 
	\begin{cases}
		1 - \frac{1}{2}\abs{r} - \abs{r}^2 + \frac{1}{2}\abs{r}^3	&	0\leq\abs{r}\leq1,	\\
		1 - \frac{11}{6}\abs{r} + \abs{r}^2 - \frac{1}{6}\abs{r}^3	&	1\leq\abs{r}\leq2,	\\
		0															&	2\leq\abs{r},
	\end{cases}
	\label{eq:base_delta}
\end{equation}
which has $\mathcal{O}(\Delta x^4)$ interpolation accuracy, $C^0$ regularity, and a support of four gridpoints \cite{Tornberg2004,Bringley2008}.  Define:
\begin{equation}
	\tilde\delta = \delta_{IB_4} * \delta_{IB_4} * \delta_{IB_4} * \delta_{IB_4},
\end{equation}
where $*$ denotes convolution.  It is clear that $\tilde\delta\in C^3$, its support is $16$ gridpoints, and it is a simple exercise to show that convolution preserves interpolation accuracy.  The formula for $\tilde\delta$ is given in \Cref{appendix:delta}.

\begin{remark}
	While our choice of $\delta$-function regularization limits the spatial discretization in this paper to fourth-order accuracy, this is not a fundamental limitation of the method.  One way to generalize to higher orders would be to construct smoother and more accurate regularizations of the $\delta$-function similar to the one that we construct, which maintain finite support for computational efficiency.  An alternative approach would be to use globally supported regularizations of the $\delta$-function that are $C^\infty$ (i.e. $\textnormal{sinc}$ or Gaussian functions) and make use of fast sinc transforms or fast Gaussian transforms in order to rapidly apply the spread and interpolation operators \cite{Greengard1991,Greengard2006}.
\end{remark}

\begin{remark}
	Our construction of $\tilde\delta$ is likely not optimal, in that there probably exist $C^3$ $\delta$-function regularizations accurate to $\mathcal{O}(\Delta x^4)$ with {\color{review_color}a} support {\color{review_color}of} less than 16 gridpoints (e.g. a $\delta$-function with $C^3$ regularity and accuracy of $\mathcal{O}(\Delta x^2)$ with support of only six gridpoints is defined in \cite{Bao}).  The construction of such a $\delta$-function would be worthwhile, allowing for coarser discretizations of problems with closely spaced boundaries and faster application of the operators $S$, $S^*$, $T_k$, and $T_k^*$.
\end{remark}

\subsection{Choice of extension operator $\mathcal{H}^k$}
\label{section:numerics:extension}

\begin{figure}
	\centering
	\hspace*{\fill}
	\begin{subfigure}[b]{0.4\textwidth}
		\centering
		\includegraphics[width=\textwidth]{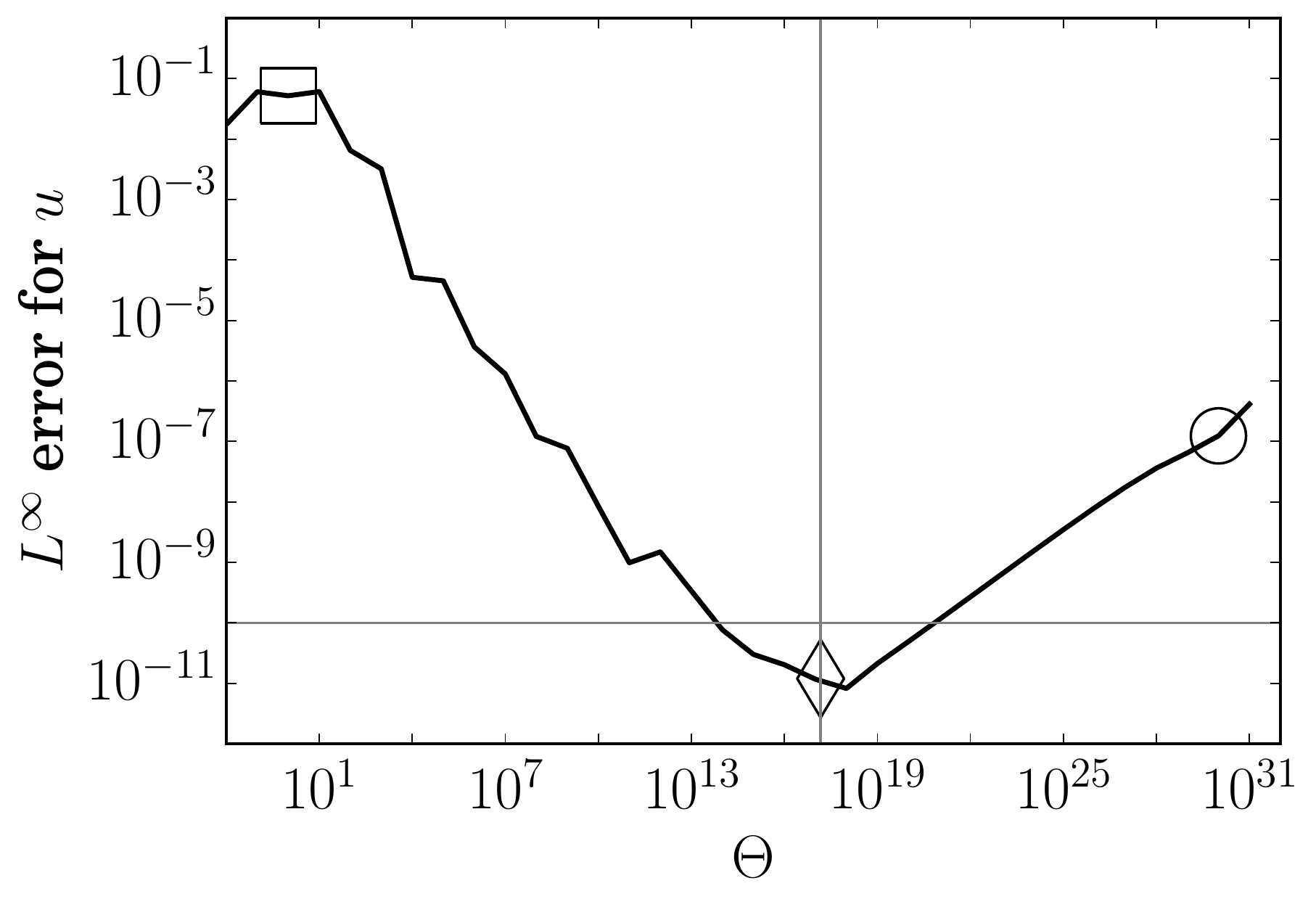}
		\subcaption{$L^\infty$ error as a function of $\Theta$}
		\label{fig:Theta_Example:plot}
	\end{subfigure}
	\hfill
	\begin{subfigure}[b]{0.4\textwidth}
	\centering
		\includegraphics[width=\textwidth]{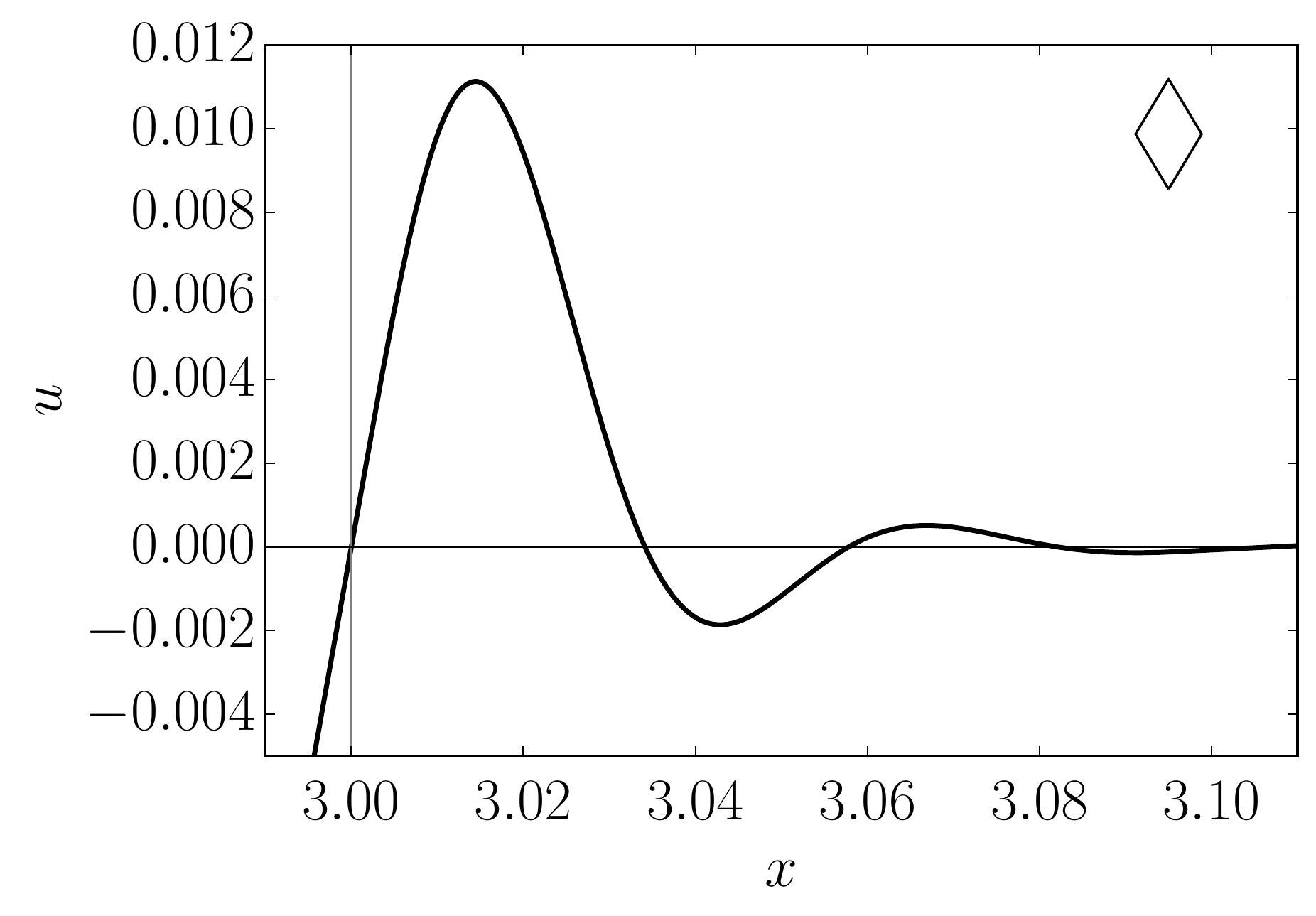}
		\subcaption{$\Theta=\Theta^*$}
		\label{fig:Theta_Example:ngood}
	\end{subfigure}
	\hspace*{\fill}
	\\
	\centering
	\hspace*{\fill}
	\begin{subfigure}[b]{0.4\textwidth}
		\includegraphics[width=\textwidth]{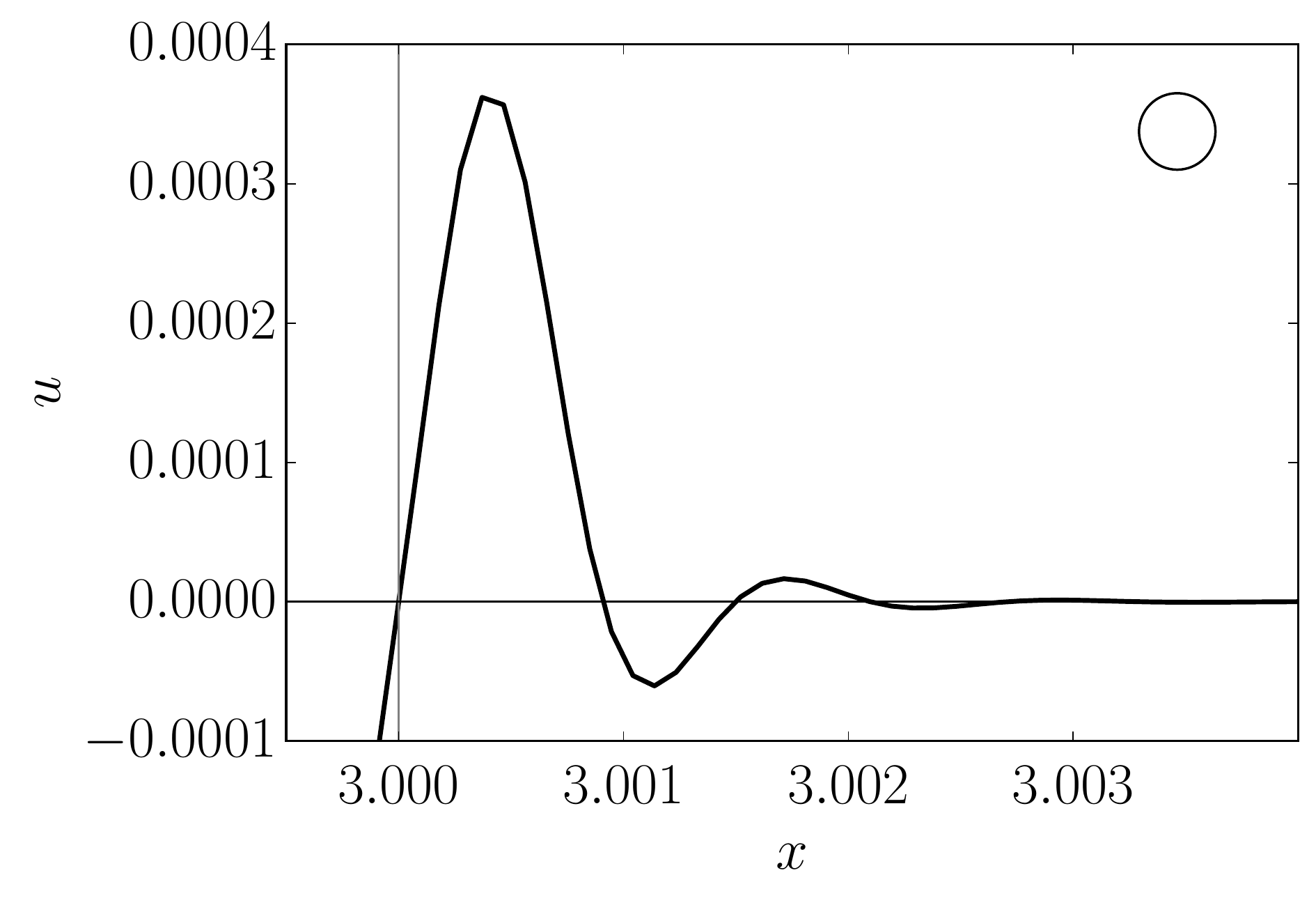}
		\subcaption{$\Theta=1\times10^{30}$}
		\label{fig:Theta_Example:nbig}
	\end{subfigure}
	\hfill
	\begin{subfigure}[b]{0.4\textwidth}
		\includegraphics[width=\textwidth]{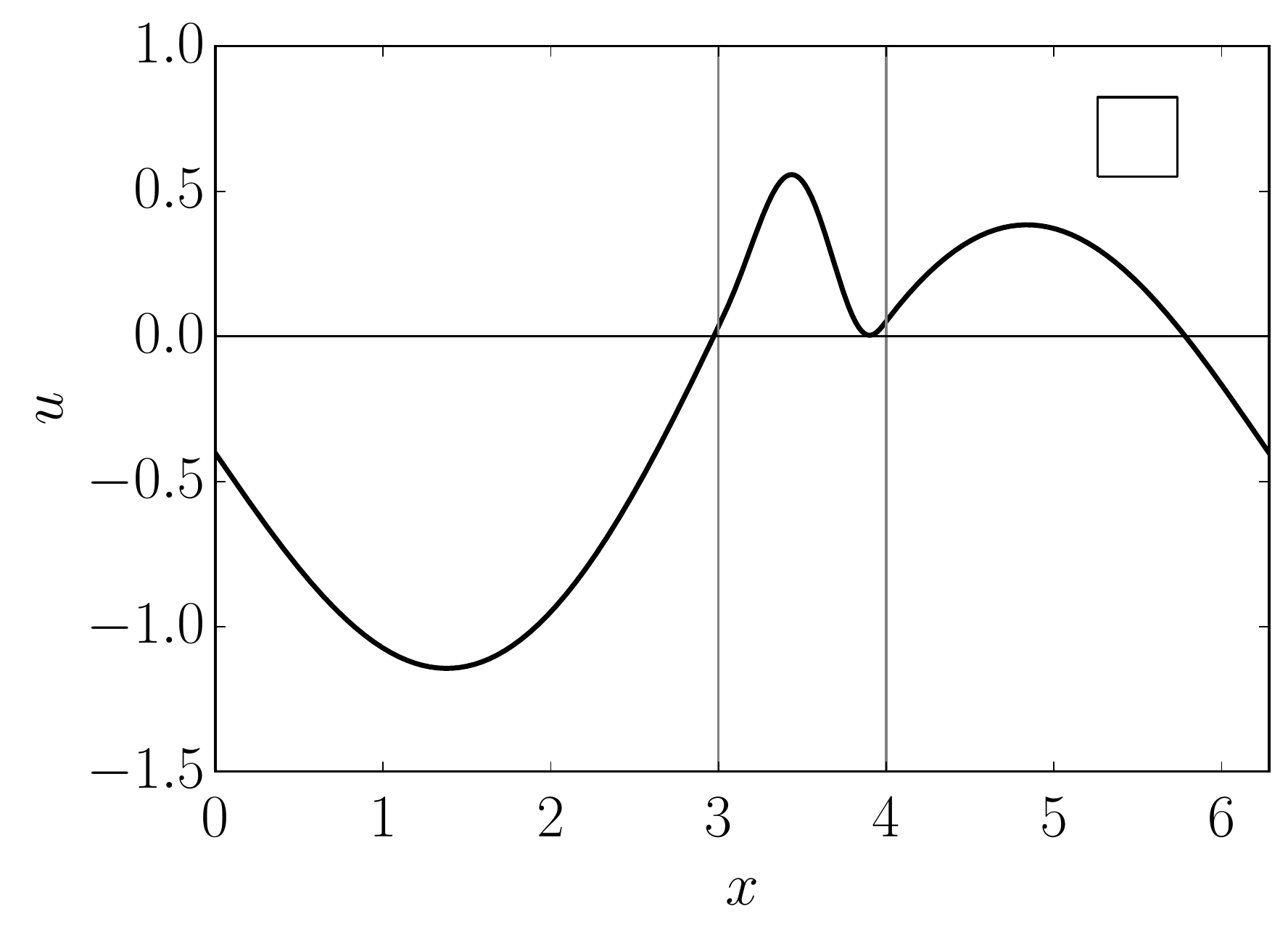}
		\subcaption{$\Theta=1$}
		\label{fig:Theta_Example:nsmall}
	\end{subfigure}
	\hspace*{\fill}
	\caption{The effect that the value $\Theta$, used to adjust the conditioning of the extension operator $\mathcal{H}^k$, has on the solution $u$ of \Cref{eq:example}.  \Cref{fig:Theta_Example:plot} shows a plot of the $L^\infty$ error against $\Theta$.  Three special points are denoted in this figure: $\Theta=\Theta^*$ as defined in \Cref{eq:Theta_definition}, denoted by ($\lozenge$); $\Theta=10^{30}$, denoted by (\Circle); and $\Theta=1$, denoted by (\Square).  The solutions corresponding to these special values of $\Theta$, along with their extensions, are shown in \Cref{fig:Theta_Example:ngood,fig:Theta_Example:nbig,fig:Theta_Example:nsmall}, respectively.  For \Cref{fig:Theta_Example:ngood,fig:Theta_Example:nbig}, only a zoom of the area near $x=3$ is shown.  See \Cref{section:numerics:extension} for a detailed discussion.}
	\label{fig:Theta_Example}
\end{figure}

Due to the nullspace of the periodic polyharmonic operator $\Delta^{k+1}$, we choose the operator $\mathcal{H}^k$ used in the extension problem {\color{review_color}given by} \Cref{eq:extension} to be the invertible operator
\begin{equation}
	\label{eq:Helmholtz_definition}
	\mathcal{H}^k = \Delta^{k+1} + (-1)^{k+1}\Theta(k,n).
\end{equation}
Here $\Theta$ is a positive scalar function that depends on the smoothness of the extension ($k$) and the number of Fourier modes ($n$) used in the discretization of the problem.  The function $\Theta$ is chosen to mitigate the numerical condition number of the operator $\mathcal{H}^k$:
\begin{equation}
	\kappa = 1 + \frac{(n/2)^{2(k+1)}}{\Theta}.
\end{equation}
Minimizing the condition number $\kappa$ (by taking $\Theta$ to be large) must be balanced against the need to resolve the intrinsic length scale $L=\Theta^{-1/2(k+1)}$ introduced to the problem.  In practice, we find that taking $\Theta$ to be
\begin{equation}
	\label{eq:Theta_definition}
	\Theta^* = \textnormal{max}\left(1,\alpha\varepsilon\left(\frac{n}{2}\right)^{2(k+1)}\right),
\end{equation}
where $\varepsilon$ is machine precision ($2^{-52}$ for double precision and $2^{-112}$ for quadruple precision) and $\alpha=0.001$ produces stable and accurate solutions.  Numerical results are not particularly sensitive to the choice of $\alpha$ (see \Cref{fig:Theta_Example:plot}).  This stability is demonstrated later (in \Cref{section:poisson_test:1d}) with highly accurate solutions shown for the IBSE-$1$, IBSE-$2$, and IBSE-$3$ methods up to $n=2^{22}$.  All numerical results in this paper are produced using $\Theta^*$ to construct $\mathcal{H}^k$.

For large values of $n$ and $k$, the value of $\Theta$ can substantially impact the overall accuracy of the IBSE-$k$ algorithm.  We demonstrate this effect in \Cref{fig:Theta_Example} for $n=2^{16}$, $k=3$, and solutions to \Cref{eq:example}.  In \Cref{fig:Theta_Example:plot}, we plot the $L^\infty$ error in the solution $u$ as a function of $\Theta$ across 32 orders of magnitude.  The minimum error observed is $7.63\times10^{-12}$, at $\Theta=10^{17}$, which is within an order of magnitude of $\Theta^*$ (denoted by a vertical grey line).  Although there is some sensitivity to $\Theta$ near the minimum, it is small: the error is bounded by $10^{-10}$ over 6 orders of magnitude of $\Theta$.  In \Cref{fig:Theta_Example:ngood}, we show a zoom of the extension near the domain boundary at $x=3$ for the solution produced using $\Theta^*$, which has an error of $1.12\times 10^{-11}$.  The extension decays rapidly to $0$, but is well resolved by the discretization: there are $150$ gridpoints between the boundary at $x=3$ and the peak at $x\approx 3.014$.  In \Cref{fig:Theta_Example:nbig}, we show a zoom of the extension produced when $\Theta=10^{30}$, which gives a solution with an error of $1.23\times10^{-7}$, four orders of magnitude worse than the error produced when $\Theta=\Theta^*$.  Although taking $\Theta$ this large leads to a very well-conditioned operator $\mathcal{H}^k$, the extension decays \emph{very} rapidly to $0$, and the discretization is no longer able to fully resolve the length scales: there are only \emph{three} gridpoints between the boundary at $x=3$ and the peak at $x\approx3.0004$.  Finally, in \Cref{fig:Theta_Example:nsmall}, we use $\Theta=1$.  Although there are no fine length-scales to be resolved, ill-conditioning in the operator $\mathcal{H}^k$ leads to an error of $5.2\times 10^{-2}$.

\subsection{Inversion of the IBSE-$k$ system \eqref{eq:the_system}}
\label{section:numerics:inversion}

In this section, we consider solving \Cref{eq:the_system} for an \emph{invertible} elliptic operator $\mathcal{L}$ in place of $\Delta$.  The case where $\mathcal{L}$ is the periodic Laplacian is complicated by the nullspace of $\Delta$: even though \Cref{eq:the_system} is invertible, the method for inversion that we outline in this section is not directly applicable; details for the resolution of this problem are provided in \Cref{appendix:poisson_inversion}.  The system of equations is
\begin{equation}
	\label{eq:the_system_2}
	\left(
	\begin{array}{cc|cc}
		\mathcal{L}	&	-\chi_E\mathcal{L}	&		&	S	\\
				&	\mathcal{H}^k	&	T_k	&		\\
		\hline
		-T_k^*	&	T_k^*			&		&		\\
		S^*		&					&		&		
	\end{array}
	\right)
	\begin{pmatrix}
		u	\\	\xi	\\	F	\\	G
	\end{pmatrix}
	=
	\begin{pmatrix}
		\chi_\Omega	f	\\	0	\\	0	\\	g
	\end{pmatrix}.
\end{equation}
This system is large: in two spatial dimensions, the number of equations in the system is $2n^2+(k+1)n_\text{bdy}$.  By block Gaussian elimination, we can find a \emph{Schur-complement} for this matrix, which we label $SC$:
\begin{equation}
	SC =
	\begin{pmatrix}
		-T_{\color{review_color}k}^*	&	T_{\color{review_color}k}^*	\\
		S^*
	\end{pmatrix}
	\begin{pmatrix}
		\mathcal{L}	&	-\chi_E\mathcal{L}	\\
					&	\mathcal{H}^k
	\end{pmatrix}^{-1}
	\begin{pmatrix}
				&	S	\\
		T_{\color{review_color}k}
	\end{pmatrix},
	\label{eq:the_schur_complement}
\end{equation}
along with an associated system of equations for the Lagrange multipliers $F$ and $G$:
\begin{equation}
	\bigg(
		\quad SC \quad
	\bigg)
	\begin{pmatrix}
		F	\\	G
	\end{pmatrix}
	=
	\begin{pmatrix}
		S^*\mathcal{L}^{-1}\chi_\Omega f - g	\\
		T_k^*\mathcal{L}^{-1}\chi_\Omega f
	\end{pmatrix}.
	\label{eq:schur_system}
\end{equation}
The size of $SC$ is comparatively small, only $(k+1)n_\text{bdy}$ square.  This equation is the key to the efficiency of the algorithm, as it allows the Lagrange multipliers $F$ and $G$ to be computed rapidly without first solving for $u$ or $\xi$.  Once $F$ and $G$ are determined, we may solve for $\xi$:
\begin{equation}
	\label{eq:get_xi}
	\xi = -(\mathcal{H}^k)^{-1}T_kF,
\end{equation}
and once $\xi$ is known, $u$ may be computed as
\begin{equation}
	\label{eq:get_u}
	u = \mathcal{L}^{-1}\left(\chi_Df + \chi_E\mathcal{L}\xi - SG\right).
\end{equation}
The computation of $\xi$ and $u$ requires only a series of fast operations: FFTs, multiplies, and applications of the spread operators.

Rapid inversion of the system of equations for $F$ and $G$ given in \Cref{eq:schur_system} is not a trivial task.  Because we have restricted to problems set on stationary domains and the size of $SC$ is small, it is feasible to proceed using dense linear algebra.  We form $SC$ by repeatedly applying it to basis vectors.  This operation is expensive: for two dimensions it is $\mathcal{O}(N^{3/2}\log N)$ in the total number of unknowns $N=n^2$.  The Schur-complement is then factored by the LU algorithm provided by LAPACK \cite{laug}.  Once this factorization is computed \Cref{eq:schur_system} can be solved rapidly.  This Schur-complement depends only on the domain and the discretization, so the LU-decomposition can be reused to solve multiple problems on the same domain or in each timestep when solving time-dependent problems.  See \Cref{section:numerics:complexity} for a discussion of the computational complexity of the method and \Cref{table:heat_equation_times} for the actual numerical cost of the method applied to a test problem.

\subsection{Summary of algorithm and numerical complexity}
\label{section:numerics:complexity}
{\color{review_color}
For a given geometry, choice of elliptic operator $\mathcal{L}$, and discretization size $n$, the implementation of the IBSE-$k$ method proceeds as follows.
\begin{enumerate}
	\item Setup
	\begin{enumerate}
		\item The boundary $\Gamma$ is discretized as described in \Cref{section:numerics:operators}, and stencils for the evaluation of $\tilde\delta$ and its first $k$ normal derivatives are computed at all nodes of the discrete boundary in $\tilde\Gamma$ to allow for the rapid application of the $S$, $S^*$, $T_k$, and $T_k^*$ operators.
		\item The Schur-complement matrix (\Cref{eq:the_schur_complement} for invertible $\mathcal{L}$, and \Cref{eq:augmented_SC} for $\mathcal{L}=\Delta$ or other non-invertible $\mathcal{L}$) is formed columnwise, by repeatedly applying it to basis vectors.  For invertible $\mathcal{L}$, this task may be completed as follows (for non-invertible $\mathcal{L}$, the process is nearly identical):
		\begin{enumerate}
			\item One element of either $F$ or $G$ is set to $1$, all other elements are set to $0$.
			\item These forces are spread to the Eulerian grid, by the application of $T_k$ and $S$.
			\item $\xi$, and then $u$, are computed by applying \Cref{eq:get_xi} followed by \Cref{eq:get_u}.
			\item The quantities $S^*u$ and $-T_k^* u + T_k^*\xi$ are computed, and placed in the appropriate rows and column of the Schur-complement.
		\end{enumerate}
		Note that item (ii) corresponds to application of the right-most matrix in \Cref{eq:the_schur_complement}, item (iii) corresponds to inverting the middle matrix in \Cref{eq:the_schur_complement}, while item (iv) corresponds to the application of the left-most matrix in \Cref{eq:the_schur_complement}.  We remark that although this operation is expensive (see \Cref{table:heat_equation_times}), it is trivially parallelizable, as the task of applying the Schur-complement to each individual basis vector is independent.
		\item The LU decomposition of the Schur-complement is computed.
		\item Optionally, these items may be saved, bypassing the expensive computation and factorization of the Schur-complement in future uses of the IBSE-$k$ method when applied to the same geometry, $\mathcal{L}$, and discretization size $n$.
	\end{enumerate}
	\item Solve
	\begin{enumerate}
		\item The right-hand side of \Cref{eq:schur_system} (or equivalently, \Cref{eq:augmented_SC} when $\mathcal{L}$ is not invertible) is computed.
		\item Using the pre-computed LU decomposition of the Schur complement matrix, the Lagrange multipliers $F$ and $G$ are found by solving \Cref{eq:schur_system} (or equivalently, \Cref{eq:augmented_SC} when $\mathcal{L}$ is not invertible).
		\item The smooth extension $\xi$, and finally $u$, are computed via \Cref{eq:get_xi,eq:get_u} (or equivalently, \Cref{eq:get_u_noninvertible} when $\mathcal{L}$ is not invertible)
	\end{enumerate}
\end{enumerate}
In items 1(b), 2(a), and 2(c), application of the FFT and IFFT is used to move variables between frequency and physical space, while differential operators are applied in frequency space and multiplications with characteristic functions are carried out in physical space.
}

The computational work for this method can be broken into two main tasks: ({\small\emph{i}}) {\color{review_color}the setup items, listed above, which are dominated by} the formation and factorization of the Schur-complement defined in \Cref{eq:the_schur_complement} and ({\small\emph{ii}}) the production of one solution $u$ given the Schur-complement, its factorization, the body forcing $f$, and the boundary conditions $g$.  We summarize the scaling of the algorithm in \Cref{table:cost}; see \Cref{table:heat_equation_times} in \Cref{section:heat:easy} for the actual numerical cost of the method applied to a test problem.  In two dimensions, the cost of a solve scales like the cost of an FFT; in three dimensions the cost is slightly worse: $\mathcal{O}(N^{4/3})$ in the total number of unknowns $N=n^3$.  The higher cost of the algorithm in three-dimensions is due to the cost of factoring and inverting the Schur-complement, however, these matrices are highly structured.  It is quite possible that the cost of inversion could be reduced with an appropriately preconditioned iterative method, by directly exploiting the structure of the matrix, e.g. with an inverse Multipole Method \cite{Ambikasaran2014}, or by indirectly exploiting that structure using an algorithm like HODLR \cite{Ambikasaran} or Algebraic Multigrid \cite{ruge1987algebraic}.

\begin{table}
	\centering
	\begin{tabular}{cccc}
		\toprule
		Dimension	&	Cost of FFT	&	Cost of {\color{review_color}Setup}	&	Cost of Solve	\\
		\midrule
		2	&	$\mathcal{O}(N\ln N)$	&	$\mathcal{O}(N^{3/2}\ln N)$	&	$\mathcal{O}(N\ln N)$	\\
		3	&	$\mathcal{O}(N\ln N)$	&	$\mathcal{O}(N^2)$	&	$\mathcal{O}(N^{4/3})$	\\
		\bottomrule
	\end{tabular}
	\caption{Numerical cost of the IBSE-$k$ algorithm as a function of $N=n^d$ for dimensions $d=2$ and $d=3$, along with the cost of an FFT.  The `Cost of Solve' refers to the cost of solving \Cref{eq:the_system} once, which is equivalent to the cost of solving one elliptic problem or taking one implicit timestep of a parabolic problem.}
	\label{table:cost}
\end{table}



\section{Results: Poisson equation}
\label{section:poisson_test}

\subsection{One-dimensional test problem}
\label{section:poisson_test:1d}

\begin{figure}
	\centering
	\includegraphics[width=0.6\textwidth]{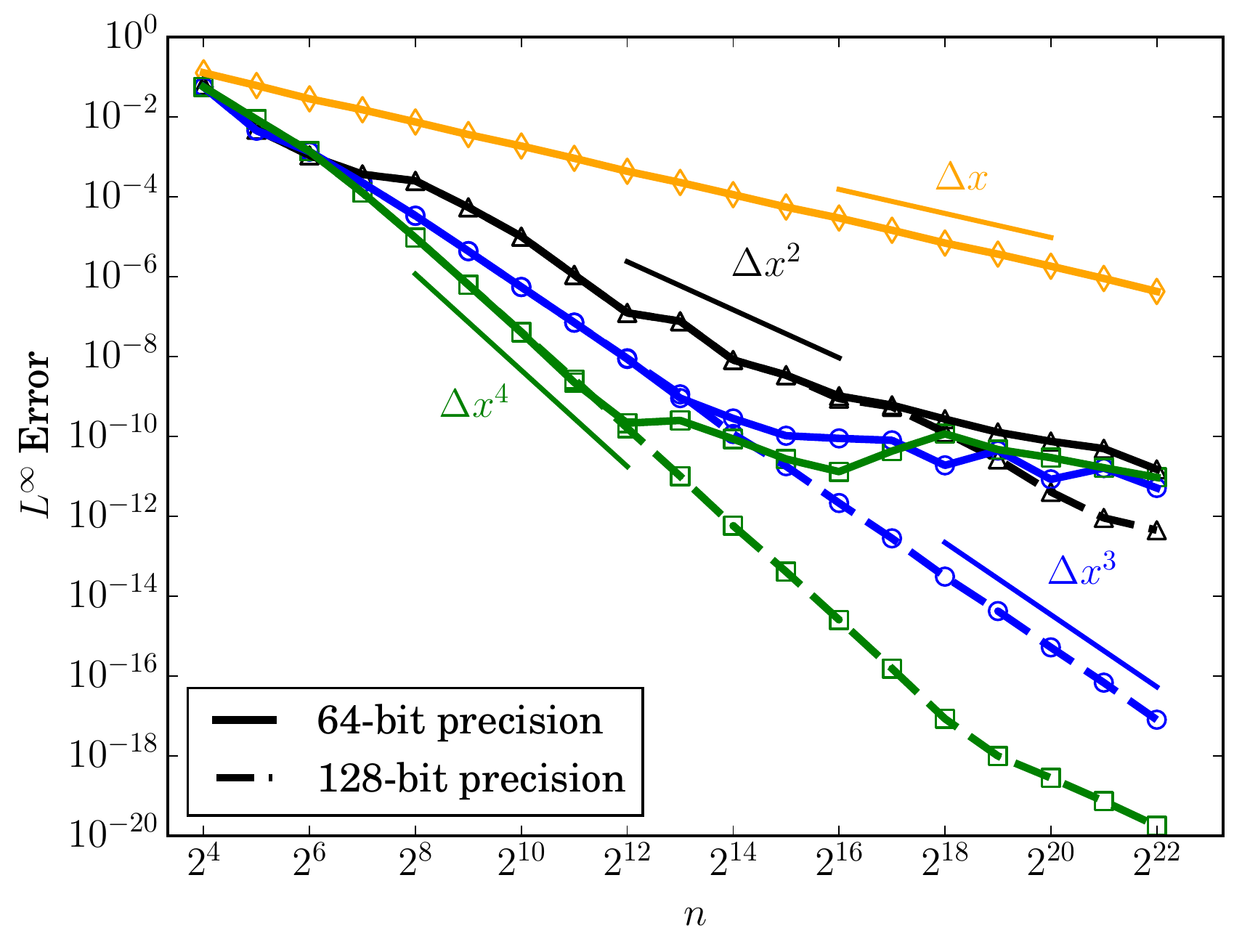}
	\caption{$L^\infty$ error for solutions to \Cref{eq:example} produced with IB ({\color{orange} $\lozenge$}), IBSE-$1$ ({\color{black} \TriangleUp}), IBSE-$2$ ({\color{py_blue} \Circle}), and IBSE-$3$ ({\color{py_green} \Square}), demonstrating $\mathcal{O}(\Delta x)$, $\mathcal{O}(\Delta x^2)$, $\mathcal{O}(\Delta x^3)$, and $\mathcal{O}(\Delta x^4)$ convergence, respectively.  Solutions computed in double precision (64-bit) are shown with solid lines; solutions computed with quadruple precision (128-bit) are shown with dashed lines.}
	\label{fig:1D_refinement}
\end{figure}

To demonstrate both the power and the limitations of the IBSE method, we compute the solution to the one-dimensional example defined in \Cref{eq:example}.  Solutions are computed on grids ranging from $n=2^4$ to $2^{22}$ points, in both double- and quadruple-precision arithmetic, for the IB, IBSE-$1$, IBSE-$2$, and IBSE-$3$ methods.  \Cref{fig:1D_refinement} shows a refinement study demonstrating the expected first-, second-, third-, and fourth-order accuracy in $L^\infty$ for the IB, IBSE-$1$, IBSE-$2$, and IBSE-$3$ methods, respectively.  The IB method achieves an error of $4.23\times10^{-7}$ with $n=2^{22}$ grid points.  Imposing additional smoothness on the solution allows this error to be matched by the IBSE-$1$ method at $n=4096$, by the IBSE-$2$ method at $n=1024$, and by the IBSE-$3$ method at $n=512$.  For the IBSE-$3$ method, this is a factor of nearly 10000 less gridpoints; equivalent to a factor of 100 million less gridpoints for two-dimensional problems.  In practice, obtaining solutions accurate to six digits with the traditional IB method is often impractical; with the IBSE method this kind of accuracy can be achieved on reasonably sized grids.

In double precision, all of the methods achieve the expected convergence rate up to some value of $n$.  The double-precision solutions begin to diverge from the quadruple-precision solutions at $n=2^{17}$ for IBSE-$1$, at $n=2^{16}$ for IBSE-$2$, and at $n=2^{12}$ for IBSE-$3$.  Using $\Theta=\Theta^*$ (see \Cref{section:numerics:extension}) to control the condition number of the extension operator $\mathcal{H}^k$ provides sufficient numerical stability at all values of $n$ to enable the computation of solutions accurate to 10-11 digits.  In quadruple precision, all methods exhibit the expected order of convergence across a wide range of values of $n$, although the fourth order method begins to fall short of the expected path of convergence for extremely fine grids ($n\geq2^{20}$).

{\color{review_color}
\begin{remark}
	We note that this test problem, as well as the test problems in \Cref{section:heat:easy,section:nonlinear:burgers,section:nonlinear:fhn}, are \emph{exterior} problems set on periodic domains.  Periodicity of solutions in the physical domain $\Omega$ is enforced naturally in these problems through the use of the Fourier basis.
\end{remark}
}

\subsection{Two-dimensional test problem: solution inside a circle}
\label{section:poisson_test:circle}
\begin{figure}
	\centering
	\includegraphics[width=0.6\textwidth]{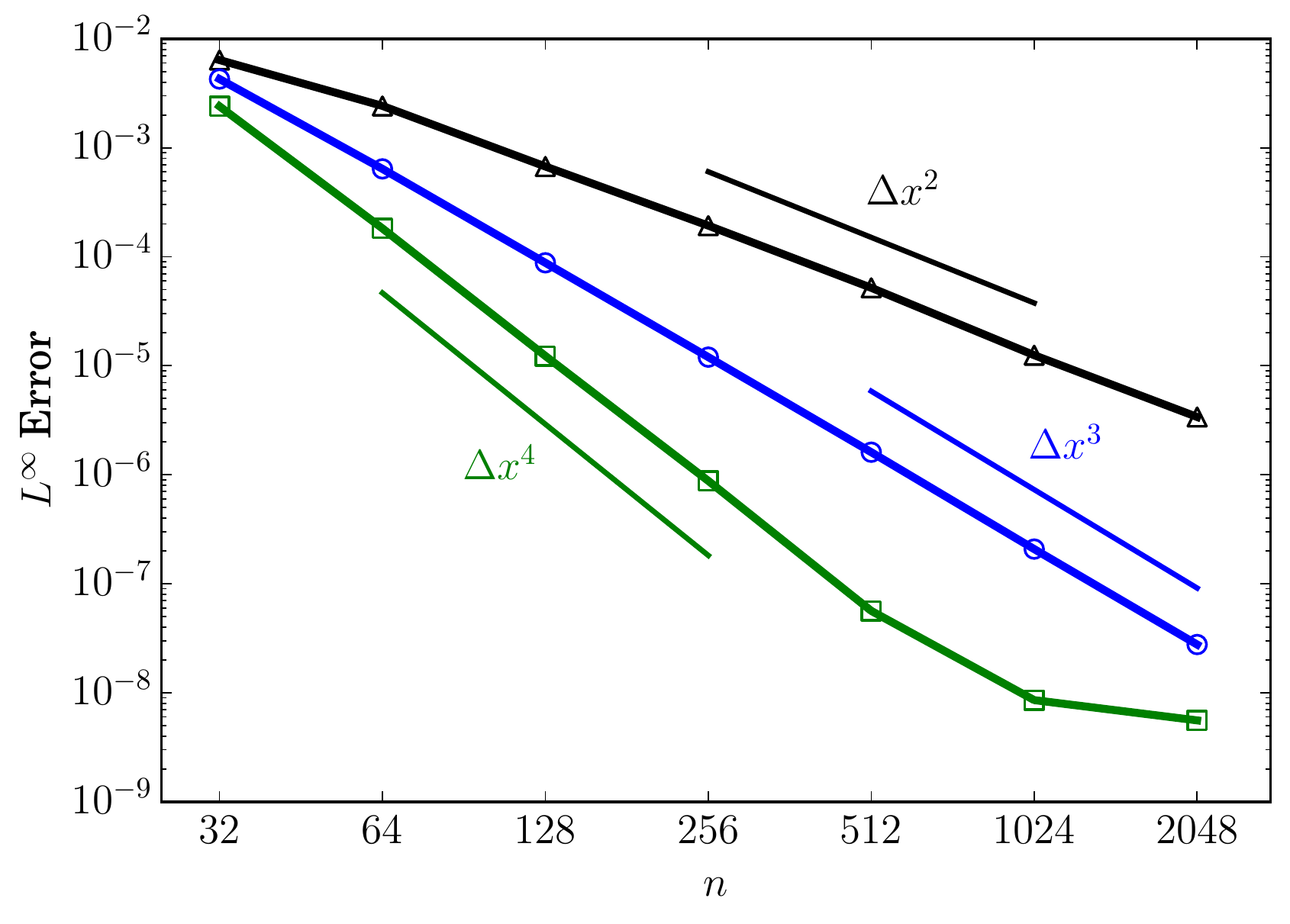}
	\caption{$L^\infty$ error for solutions to \Cref{eq:2d_unman} produced with IBSE-$1$ ({\color{black} \TriangleUp}), IBSE-$2$ ({\color{py_blue} \Circle}), and IBSE-$3$ ({\color{py_green} \Square}), demonstrating $\mathcal{O}(\Delta x^2)$, $\mathcal{O}(\Delta x^3)$, and $\mathcal{O}(\Delta x^4)$ convergence, respectively.}
	\label{fig:2D_refinement_unman}
\end{figure}

Let $\Omega=B_2((\pi,\pi))$ be the circle of radius $2$ centered at $(\pi,\pi)$, and identify $C=\mathbb{T}^2$ with the square $[0,2\pi]\times[0,2\pi]$.  We will solve the Dirichlet problem:
\begin{subequations}
	\label{eq:2d_unman}
	\begin{align}
		\Delta u	&=	-4	&	&\text{in }\Omega,	\\
		u			&=	0	&	&\text{on }\Gamma=\partial\Omega,
	\end{align}
\end{subequations}
which has an analytic solution $u_a$ given by
\begin{equation}
	u_a = 4 - (x-\pi)^2 - (y-\pi)^2.
\end{equation}
We solve this problem in double-precision arithmetic with the IBSE-$1$, IBSE-$2$, and IBSE-$3$ methods for $n=2^4$ to $n=2^{11}$.  In \Cref{fig:2D_refinement_unman}, we show $L^\infty$ error as a function of $n$, indicating second-, third-, and fourth-order accuracy for IBSE-$1$, IBSE-$2$, and IBSE-$3$ respectively.

\begin{figure}
	\centering
	\hspace*{\fill}
	\begin{subfigure}[b]{0.66\textwidth}
		\includegraphics[width=\textwidth]{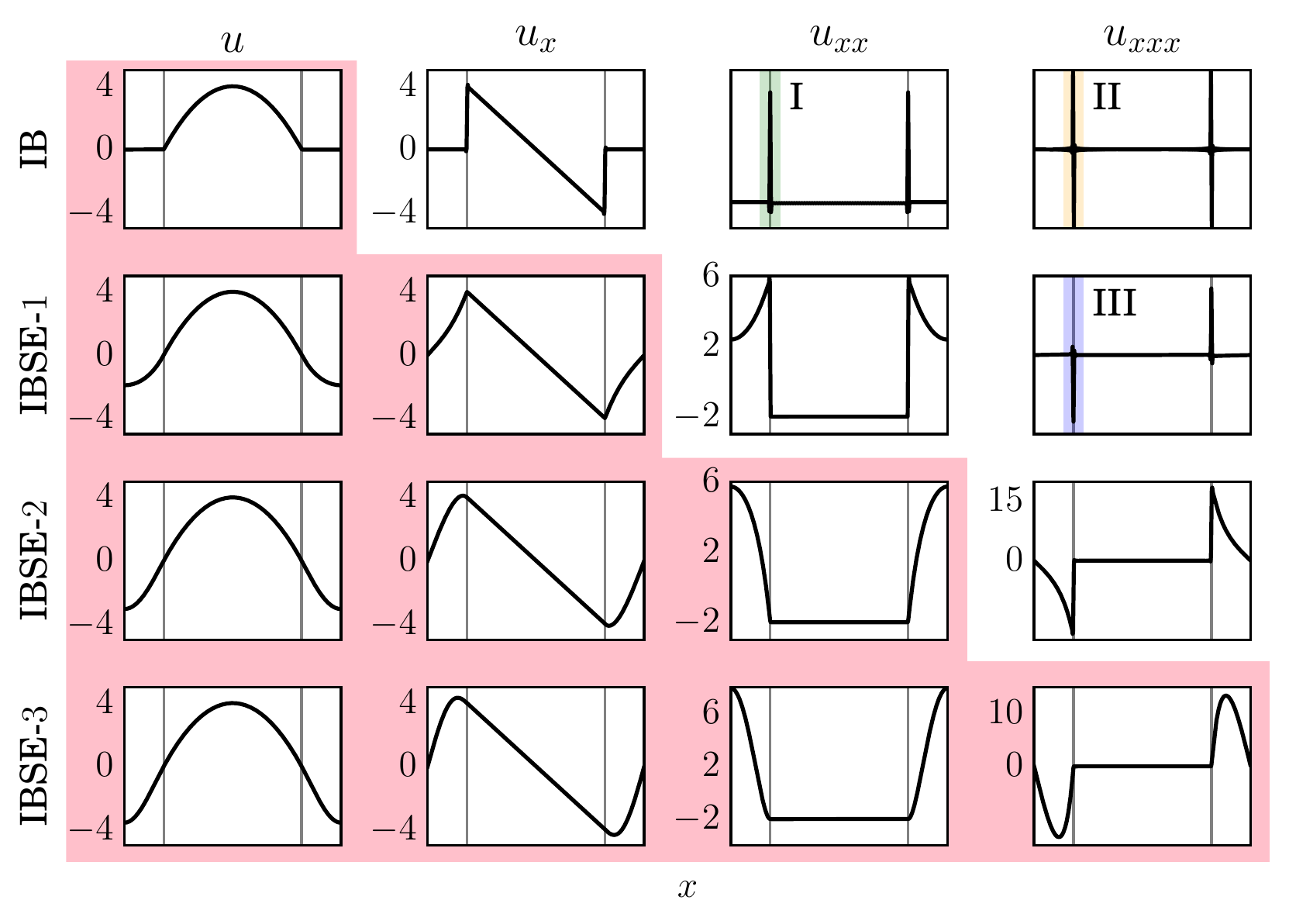}
		\caption{}
		\label{fig:2d_smooth_example:a}
	\end{subfigure}
	\hfill
	\begin{subfigure}[b]{0.32\textwidth}
		\includegraphics[width=\textwidth]{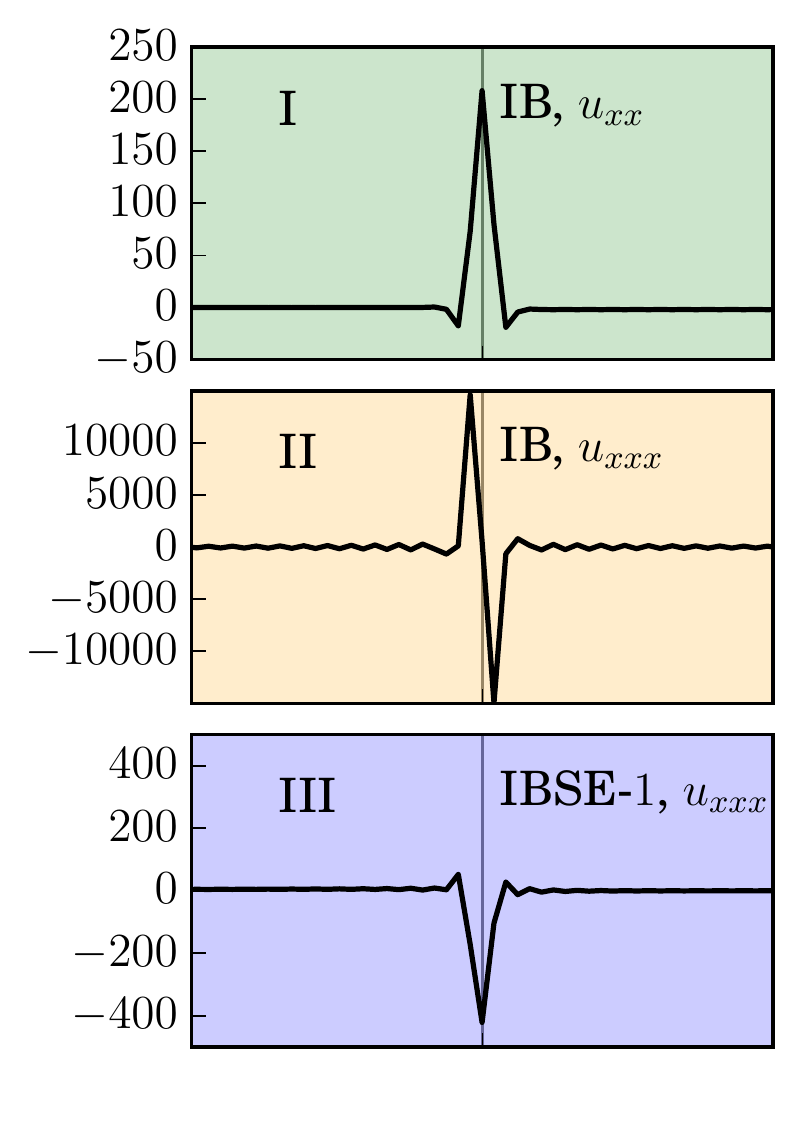}
		\caption{}
		\label{fig:2d_smooth_example:b}
	\end{subfigure}
	\hspace*{\fill}
	\caption{A demonstration of the smoothness of the solution $u$ to \Cref{eq:2d_unman} produced by the IB, IBSE-$1$, IBSE-$2$, and IBSE-$3$ methods.  \Cref{fig:2d_smooth_example:a} shows the solution $u$ and the solution's first, second, and third $x$ derivatives.  All plots show only the slice $y=\pi$ from $x=0$ to $x=2\pi$ and vertical gray lines denote the intersection of that slice with the boundary $\Gamma$ at $x=\pi-2$ and $x=\pi+2$.  The functions in plots with shading behind them (to the lower left) are \emph{at least continuous}, all other functions shown are \emph{discontinuous}.  Solutions produced using the IBSE-$k$ method exhibit discontinuities in the $(k+1)^\text{st}$ normal derivative across $\Gamma$; all lower-order derivatives are continuous.  Scales of the three plots in the upper right have been omitted.  For these plots, zooms of the region near the boundary at $x=\pi-2$ are shown in \Cref{fig:2d_smooth_example:b}.  The $x$-axis in these plots runs from $x=\pi-2-0.3$ to $x=\pi-2+0.3$.}
	\label{fig:2d_smooth_examples_2}
\end{figure}

Using the solutions to this problem, we demonstrate the property that the IBSE method produces \emph{globally smooth} solutions.  In \Cref{fig:2d_smooth_examples_2}, we show the solution $u$ to \Cref{eq:2d_unman}, along with its first, second, and third derivatives in the $x$-direction, produced using the IB, IBSE-$1$, IBSE-$2$, and IBSE-$3$ methods with $n=2^9$.  For simplicity, all functions are shown only along the slice $y=\pi$, and the intersection of this slice with the boundary $\Gamma$ is shown as gray vertical lines.  Derivatives are computed spectrally using the FFT.  We expect solutions produced by the IBSE-$k$ method to exhibit discontinuities in the $(k+1)^\text{st}$ normal derivative across interfaces, and for all derivatives of lower order to be continuous.  We see this expectation validated in the solutions presented in \Cref{fig:2d_smooth_examples_2}: functions shown in plots with shading behind them in \Cref{fig:2d_smooth_example:a} are \emph{at least continuous}, all others are discontinuous.  This property of \emph{global smoothness of the solutions} enables the Fourier-spectral discretization to obtain high-order accuracy.

\subsection{Two-dimensional test problem with Neumann boundary conditions}
\label{section:poisson_test:neumann}

Traditional Immersed Boundary approaches are unable to provide solutions to Neumann problems since convolution-style estimation of the normal derivative is not convergent at the boundary due to the lack of regularity of the solution.  The IBSE method does not have this limitation because of the additional smoothness of the solution.  The only modification necessary to adapt \Cref{eq:the_system} to solve Neumann problems is to change the $S$ and $S^*$ operators that enforce and specify the boundary conditions.  Consider the general Neumann problem:
\begin{subequations}
	\begin{align}
		\Delta u						&=	f	&	&\text{in }\Omega,	\\
		\frac{\partial u}{\partial n}	&=	g	&	&\text{on }\Gamma.
	\end{align}
\end{subequations}
Define
\begin{subequations}
	\begin{align}
		T_{(1)}F(x) &= \int_\Gamma F_j(s)\frac{\partial\delta(x-s)}{\partial n}\,ds,	\\
		T^*_{(1)}\xi(s) &= -\int \xi(x)\frac{\partial\delta(x-s)}{\partial n}\,dx, \quad\forall s\in\Gamma.
	\end{align}
\end{subequations}
The analogue to \Cref{eq:the_system} is
\begin{equation}
	\label{eq:the_system_neumann}
	\begin{pmatrix}
		\Delta	&	-\chi_E\Delta	&		&	T_{(1)}	\\
				&	\mathcal{H}^k	&		&		\\
		-T_k^*	&	T_k^*				&		&		\\
		T_{(1)}^*		&					&		&		
	\end{pmatrix}
	\begin{pmatrix}
		u	\\	\xi	\\	F	\\	G
	\end{pmatrix}
	=
	\begin{pmatrix}
		\chi_\Omega	f	\\	0	\\	0	\\	g
	\end{pmatrix}.
\end{equation}

\begin{remark}
	It is just as simple to modify the IBSE method to allow the solution of Robin problems.  Rather than replacing the upper-right and lower-left operators in \Cref{eq:the_system} with $T_{(1)}$ and $T^*_{(1)}$, they should instead be replaced with linear combinations of $S$, $T_{(1)}$, $S^*$, and $T^*_{(1)}$.
\end{remark}

\begin{figure}
	\centering
	\includegraphics[width=0.6\textwidth]{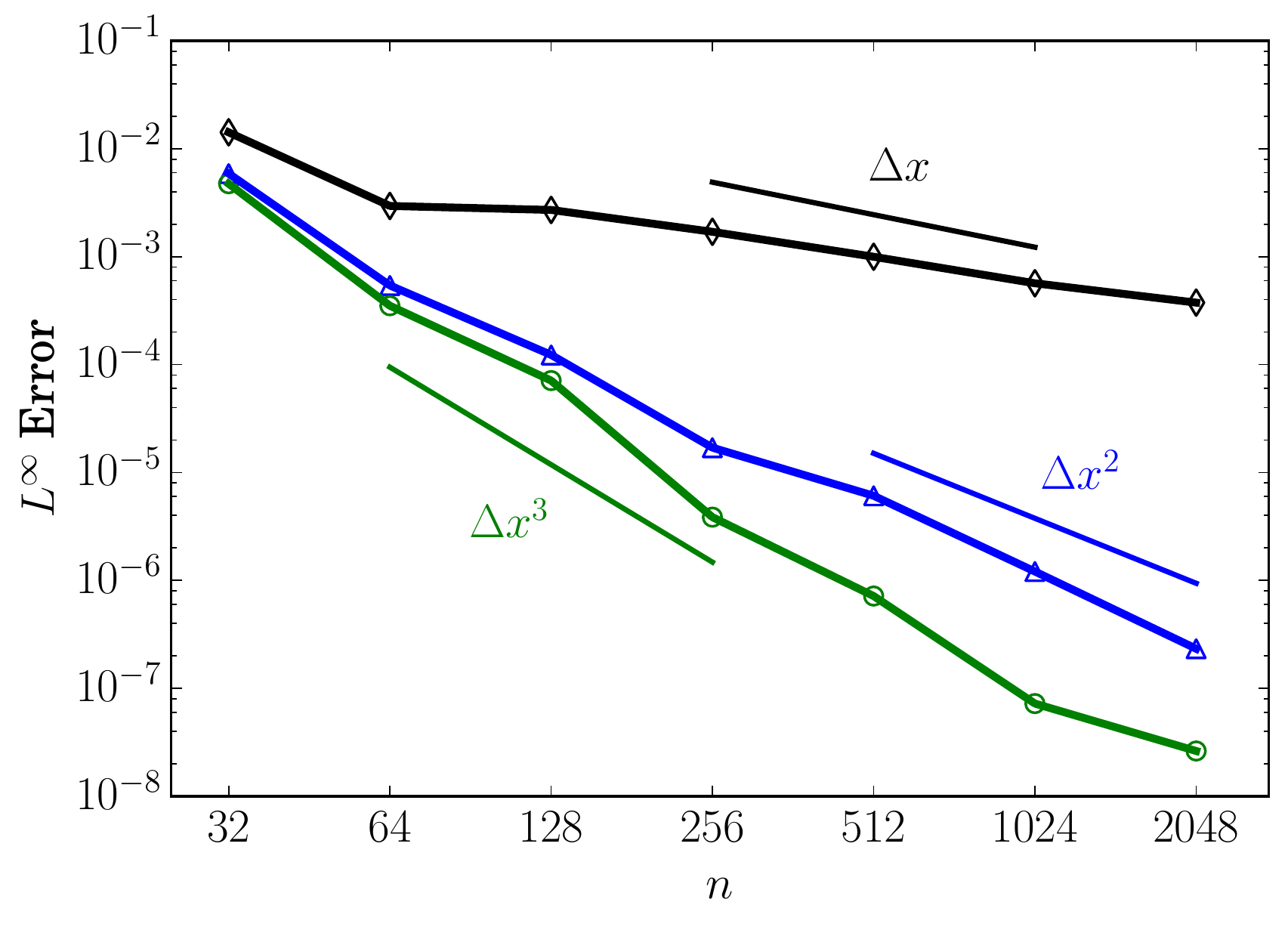}
	\caption{$L^\infty$ error for solutions to \Cref{eq:2d_neumann} produced using IBSE-$1$ ({\color{black} \TriangleUp}), IBSE-$2$ ({\color{py_blue} \Circle}), and IBSE-$3$ ({\color{py_green} \Square}) demonstrating  slightly less than first-order convergence, second-order convergence, and third-order convergence, respectively.}
	\label{fig:2D_refinement_neumann}
\end{figure}

Let $\Omega=B_1((\pi,\pi))$.  We will solve the problem
\begin{subequations}
	\label{eq:2d_neumann}
	\begin{align}
		\Delta u	&=	e^{\sin x}\left(\cos^2 x - \sin x\right) - \cos y	&	&\text{in }\Omega,	\\
		\frac{\partial u}{\partial n}	&=	\cos^2 x\, e^{\sin x} - \sin^2 y	&	&\text{on }\Gamma,
	\end{align}
\end{subequations}
which has the analytic solution
\begin{equation}
	u_a = e^{\sin x} + \cos y.
\end{equation}
We solve this problem with the IBSE-$1$, IBSE-$2$, and IBSE-$3$ methods, for which we expect first-, second-, and third-order convergence, respectively.  The reason for the lower order of convergence than in Dirichlet problems is that our discretization of $T_{(1)}^*$ is one order less accurate than $S^*$ when acting on functions of the same regularity.  Despite this, the solution $u$ produced by IBSE-$k$ still has global $C^k$ regularity.  Solutions are computed in double precision, for $n=2^4$ to $n=2^{11}$.  \Cref{fig:2D_refinement_neumann} shows $L^\infty$ error as a function of $n$.  Second- and third-order convergence is observed for the IBSE-$2$ and IBSE-$3$ methods.  The IBSE-$1$ method converges at a rate that is slightly less than first-order.

\section{Results: heat equation}
\label{section:heat}

We now consider solving the heat equation with Dirichlet boundary conditions:
\begin{subequations}
	\begin{align}
		u_t - \Delta u	&=	f(x,t)	&	&\text{in }\Omega,	\\
		u				&=	g(x,t)	&	&\text{on }\Gamma.
		\label{eq:the_heat}
	\end{align}
\end{subequations}
Since the IBSE method can directly invert elliptic equations, it can easily be adapted to provide high-order implicit-time discretizations for the heat equation.  Many implicit methods could be used; we choose BDF4 \cite{suli2003introduction}:
\begin{equation}
	\left(\mathbb{I}-\frac{12}{25}\Delta t\Delta\right)u^{n+1} = 
		\frac{1}{25}\left(12\Delta t f(t+\Delta t) + 48u^n-36u^{n-1}+16u^{n-2}-3u^{n-3}\right).
\end{equation}
Applying this scheme to \Cref{eq:the_heat} is equivalent to solving the problem:
\begin{subequations}
	\begin{align}
		\mathcal{L}u^{n+1}	&=	\tilde f		&&\text{in }\Omega,	\\
		u^{n+1}				&=	g(x,t+\Delta t)	&&\text{on }\Gamma,
	\end{align}
	\label{eq:heat_advance}
\end{subequations}
where we define:
\begin{subequations}
	\begin{align}
		\mathcal{L} &= \mathbb{I}-\frac{12}{25}\Delta t\Delta,	\\
		\tilde f		&=	\frac{1}{25}\left(12\Delta t f(t+\Delta t) + 48u^n-36u^{n-1}+16u^{n-2}-3u^{n-3}\right).
	\end{align}
\end{subequations}
This may be solved using the algorithm described in \Cref{section:numerics:inversion}.  There are several comments worth making regarding this discretization:
\begin{enumerate}
	\item The use of an implicit scheme allows for large timesteps to be taken.  For all test problems presented in this section, we choose $\Delta t$ to be:
	\begin{equation}
		\Delta t = \frac{t_\text{final}}{2}\left\lceil \frac{t_\text{final}}{\Delta x}\right\rceil^{-1},
	\end{equation}
	where $\lceil\cdot\rceil$ is the \emph{ceiling} function, so that $\Delta t\approx\Delta x / 2$.
	\item The fact that $\mathcal{L}$ depends on $\Delta t$ implies that the Schur-complement depends on $\Delta t$ as well.  This means that the Schur-complement must be reformed and refactored whenever $\Delta t$ is changed, complicating the use of a method that uses adaptive timestepping.
	\item Modifying this formulation to solve the heat equation with Neumann and Robin boundary conditions may be done in the same way as shown for the Poisson equation in \Cref{section:poisson_test:neumann}.
\end{enumerate}

\subsection{Solution in a periodic domain outside an obstacle}
\label{section:heat:easy}

\begin{figure}
	\centering
	\includegraphics[width=0.6\textwidth]{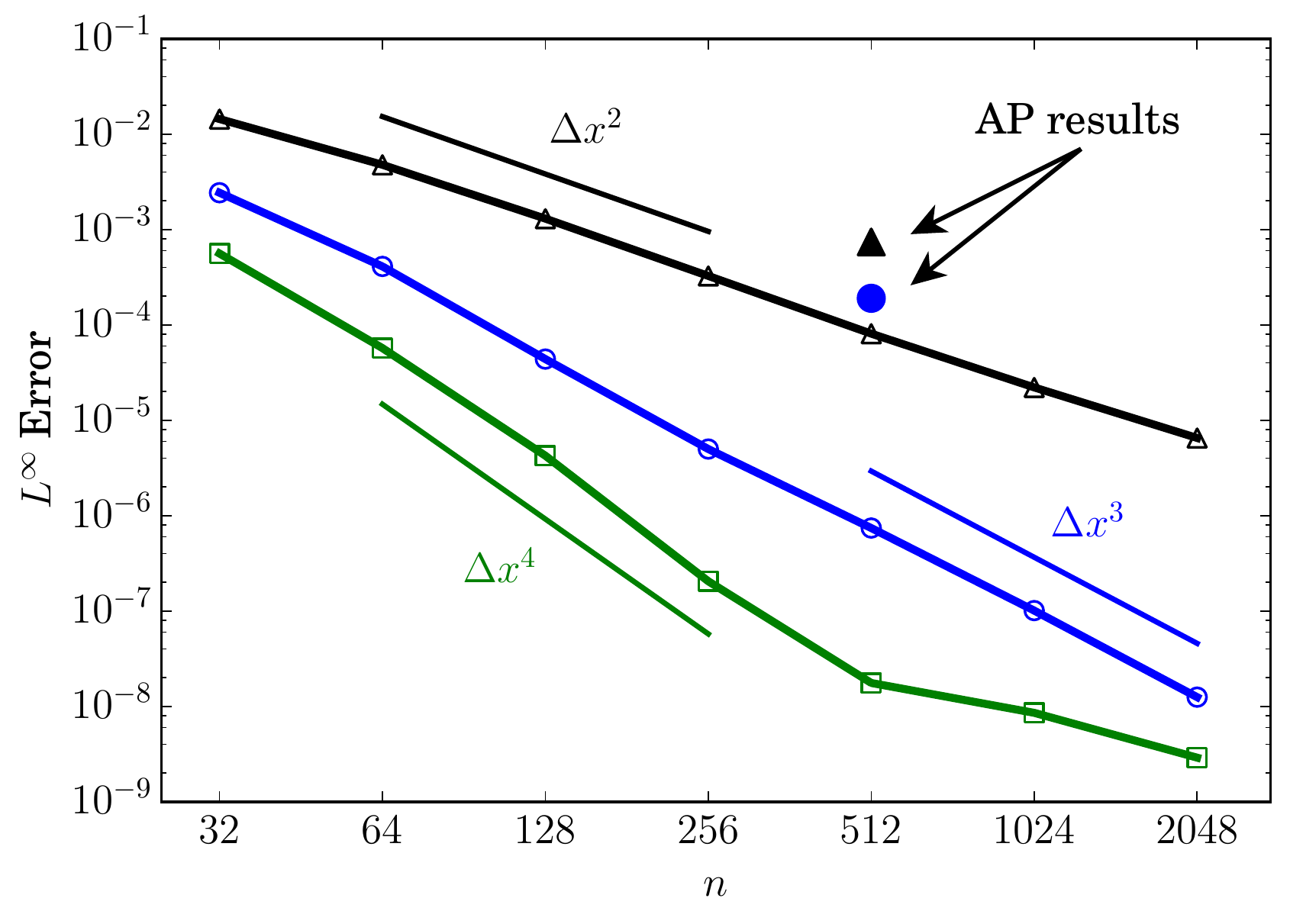}
	\caption{$L^\infty$ error for solutions to \Cref{eq:2D_test_heat} produced using IBSE-$1$ ({\color{black} \TriangleUp}), IBSE-$2$ ({\color{py_blue} \Circle}), and IBSE-$3$ ({\color{py_green} \Square}), demonstrating $\mathcal{O}(\Delta x^{2})$, $\mathcal{O}(\Delta x^3)$, and $\mathcal{O}(\Delta x^4)$ convergence, respectively.  Results from \cite{Shirokoff2013} are shown for $n=512$, for the second- (\FilledTriangleUp) and third-order({\color{py_blue} \FilledCircle}) Active Penalty methods \cite{Shirokoff2013}.}
	\label{fig:2D_refinement_heat}
\end{figure}

In this section we demonstrate high-order convergence to a heat-equation set on a periodic domain with a simple obstacle in it.  To provide direct numerical comparison with results from the Active Penalty (AP) method, the following problem is from Section 6.2 of \cite{Shirokoff2013}:
\begin{subequations}
	\label{eq:2D_test_heat}
	\begin{align}
		u_t - \Delta u	&=	
			\left[\cos y + e^{\sin x}\left(\sin x - \cos^2 x\right)\right]\cos t
			- \left(e^{\sin x} + \cos y\right)\sin t
		&	&\text{in }\Omega,	\\
		u(x,y,t=0)	&=	e^{\sin x}+\cos y		&	&\text{in }\Omega,	\\
		u		&=	(e^{\sin x}+\cos y)\cos t		&	&\text{on }\Gamma,
	\end{align}
\end{subequations}
where the physical domain is $\Omega=\mathbb{T}^2\setminus B_{1/4}(\pi,\pi)$.  \Cref{eq:2D_test_heat} has the analytic solution
\begin{equation}
	u_a=(e^{\sin x}+\cos y)\cos t.
\end{equation}
The initial condition is integrated from $t=0$ to $t=0.1$ by solving \Cref{eq:heat_advance} at each timestep.  Startup values for BDF4 are computed by evaluating the analytic solution at $t=-\Delta t$, $t=-2\Delta t$, and $t=-3\Delta t$.  In \Cref{fig:2D_refinement_heat}, $L^\infty$ errors at $t=0.1$ for solutions generated with the IBSE-$1$, IBSE-$2$, and IBSE-$3$ methods are shown, demonstrating second-, third-, and fourth-order convergence in space and time, respectively.  On this plot, we also show errors from the second- and third-order AP methods with $n=512$, the finest resolution reported in \cite{Shirokoff2013}.  Our method yields errors lower by several orders of magnitude.

In contrast to the AP method \cite{Shirokoff2013}, we are able to use an \emph{implicit} time discretization, enabling the use of large timesteps.  The ability to take large timesteps efficiently comes at the cost of a significant precomputation for the IBSE method, a cost that the AP method does not incur.  This precomputation prevents the IBSE method as described in this paper from being efficiently applied to moving boundary problems, a limitation not faced by the AP method.

In \Cref{table:heat_equation_times}, we provide the computational time required by the IBSE-$1$, IBSE-$2$, and IBSE-$3$ methods to setup the problem (dominated by the formation and factorization of \Cref{eq:the_schur_complement}), the time required to take one timestep, and the number of timesteps taken to advance the equation from $t=0$ to $t=0.1$.  Computational times are given as a multiple of the time required to take one FFT.  In our code, 8 FFTs are used per timestep regardless of $k$; the time required to take a timestep does not vary significantly across the methods.  The precomputation time increases approximately linearly in $k$, but even for the finest discretization presented, the precomputation time is not prohibitive:  for $n=2048$ and the IBSE-$3$ method, the precomputation requires about 50 minutes using serial code that has not been carefully optimized.

\setlength{\tabcolsep}{4pt}
\begin{table}
	\centering
		\begin{tabular}{clcclcclcclc}
		\toprule
		\multirow{2}{*}{$n\times n$} && \multicolumn{2}{c}{IBSE-$1$} && \multicolumn{2}{c}{IBSE-$2$} && \multicolumn{2}{c}{IBSE-$3$}	&&	\multirow{2}{*}{\parbox{2cm}{\centering timesteps to $t=0.1$}} \\
		\cmidrule{3-4} \cmidrule{6-7} \cmidrule{9-10}
					&&	prep time	&	timestep	&&	prep time	&	timestep	&&	prep time	&	timestep	&&		\\
		\midrule
		$32\times 32$		&&	674			&	39			&&	895			&	39			&&	1451		&	43			&&	2		\\
		$64\times 64$		&&	957			&	34			&&	1776		&	39			&&	2371		&	39			&&	3		\\
		$128\times 128$ 	&&	1310		&	30			&&	2207		&	31			&&	3241		&	37			&&	5		\\
		$256\times 256$		&&	1833		&	21			&&	2833		&	24			&&	3917		&	23			&&	9		\\
		$512\times 512$		&&	3478		&	21			&&	4493		&	22			&&	5782		&	20			&&	17		\\
		$1024\times 1024$	&&	5328		&	19			&&	7521		&	19			&&	10187		&	19			&&	33		\\
		$2048\times 2048$	&&	9951		&	18			&&	13707		&	18			&&	18686		&	18			&&	66		\\
		\bottomrule
		\end{tabular}
	\caption{The computational time, normalized by the time to compute an FFT, required for precomputation (labeled `prep time'), and per timestep of \Cref{eq:2D_test_heat} (labeled `timestep').  We note that even for $n=2048$, the `prep time', shown here as requiring 18686 FFTs, requires only 50 minutes of wall time using serial code that has not been carefully optimized.  See text for further discussion.  The number of timesteps taken to advance \Cref{eq:2D_test_heat} from $t=0$ to $t=0.1$ is also shown.}
	\label{table:heat_equation_times}
\end{table}
\setlength{\tabcolsep}{6pt}

\subsection{Solution inside a parametrically defined region}

\begin{figure}
	\centering
	\hspace*{\fill}
	\begin{subfigure}[b]{0.48\textwidth}
		\centering
		\includegraphics[width=\textwidth]{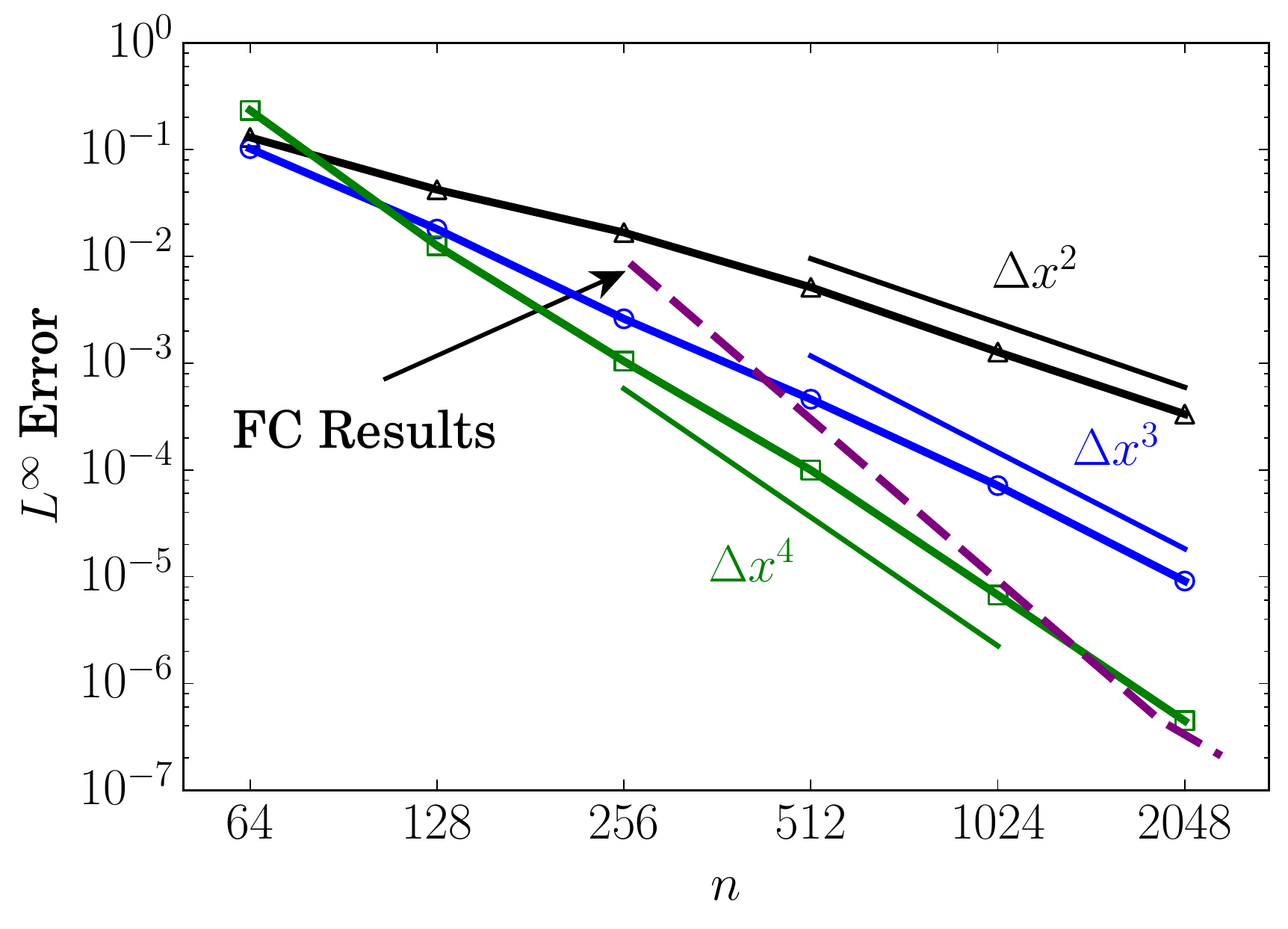}\\
		\vspace{0.05in}
		\subcaption{refinement study}
		\label{fig:2D_refinement_heat_hard:refinement}
	\end{subfigure}
	\hfill
	\begin{subfigure}[b]{0.48\textwidth}
	\centering
		\includegraphics[width=\textwidth]{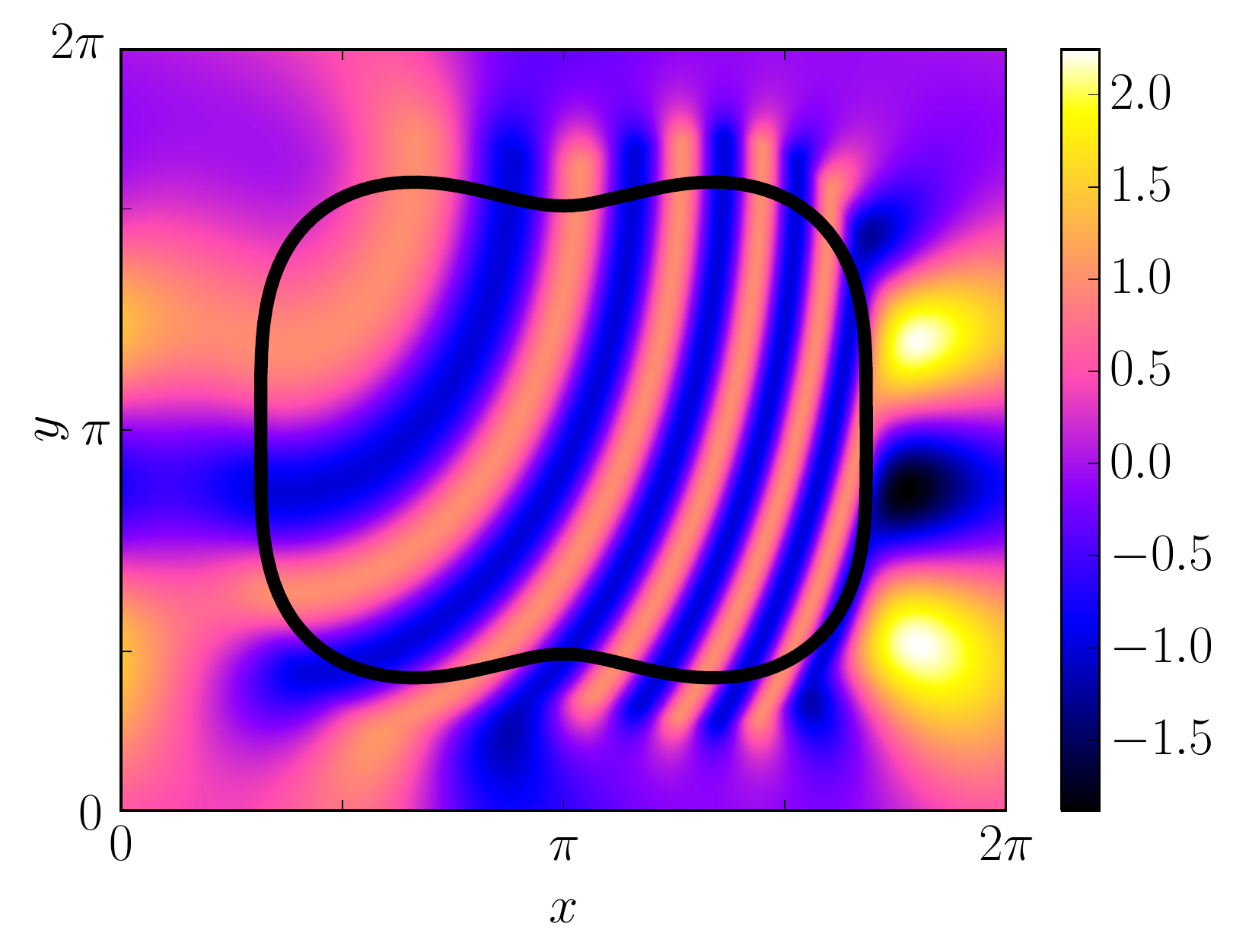}
		\subcaption{solution with extension}
		\label{fig:2D_refinement_heat_hard:solution}
	\end{subfigure}
	\hspace*{\fill}
	\caption{\Cref{fig:2D_refinement_heat_hard:refinement} shows $L^\infty$ errors for \Cref{eq:2D_test_heat_hard} at $t=1.0$ produced by the IBSE-$1$ ({\color{black} \TriangleUp}), IBSE-$2$ ({\color{py_blue} \Circle}), and IBSE-$3$ ({\color{py_green} \Square}) methods, demonstrating $\mathcal{O}(\Delta x^{2})$, $\mathcal{O}(\Delta x^3)$, and $\mathcal{O}(\Delta x^4)$ convergence, respectively.  The dashed line shows results from the Fourier Continuation method \cite{Lyon2010a}; the value of `$n$' at which we plot their results is set so that the grid-spacing $\Delta x$ is the same for the FC and IBSE methods, taking into account our rescaling of the domain.  \Cref{fig:2D_refinement_heat_hard:solution} shows the numerical solution to \Cref{eq:2D_test_heat_hard} computed with the ISBE-$1$ method and $n=2^9$ at $t=1.0$.  The boundary $\Gamma=\partial\Omega$ of the domain $\Omega$ is shown in black.  Note that the physical domain $\Omega$ is contained inside the boundary.}
	\label{fig:2D_refinement_heat_hard}
\end{figure}

In this section, we demonstrate high-order convergence to a heat equation set inside a parametrically defined domain.  This problem is from Section 6.1 of \cite{Lyon2010a}, allowing direct comparison of the IBSE method with the Fourier Continuation method \cite{Lyon2010a}.  For convenience, we rescale the domain of the problem.  We will solve
\begin{subequations}
	\label{eq:2D_test_heat_hard}
	\begin{align}
		u_t - \Delta u	&=	
			\pi\left(2-\Delta\phi\right)\cos(\pi\phi) + \pi^2\left(\phi_x^2+\phi_y^2\right)\sin(\pi\phi)	&	&\text{in }\Omega,	\\
		u(x,y,t=0)	&=	\sin(\pi\phi(x,y,t=0))		&	&\text{in }\Omega,	\\
		u		&=	\sin(\pi\phi)		&	&\text{on }\Gamma,
	\end{align}
\end{subequations}
where
\begin{equation}
	\phi(x,y,t) = 9\left(\frac{x-\pi}{4}+1\right)^2 + 4\left(\frac{y-\pi}{4}+1\right)^2 + 2t.
\end{equation}
The domain $\Omega$ is the region inside the boundary $\Gamma$ defined by the parametric equations:
\begin{subequations}
	\begin{align}
		x(\theta) &= \left(10\sin^2(2\theta) + 3\cos^3(2\theta) + 40\right)\cos\theta / 20 + \pi,	\\
		y(\theta) &= \left(10\sin^2(2\theta) + 3\cos^3(2\theta) + 40\right)\sin\theta / 20 + \pi,
	\end{align}
\end{subequations}
for $0\leq\theta<2\pi$.  \Cref{eq:2D_test_heat_hard} has the analytic solution $u=\sin(\pi\phi)$.  The initial condition is integrated from $t=0$ to $t=0.01$ by solving \Cref{eq:heat_advance} at each timestep.  Startup values for BDF4 are computed by evaluating the analytic solution at $t=-\Delta t$, $t=-2\Delta t$, and $t=-3\Delta t$.  In \Cref{fig:2D_refinement_heat_hard}, we report the maximum error over space and time, demonstrating second-, third-, and fourth-order convergence in space and time with the IBSE-$1$, IBSE-$2$, and IBSE-$3$ methods, respectively.  \Cref{fig:2D_refinement_heat_hard:solution} shows a solution to \Cref{eq:2D_test_heat_hard}, together with its $C^1$ extension, computed using IBSE-$1$ and $n=2^9$ at $t=1.0$.  Also shown in \Cref{fig:2D_refinement_heat_hard:refinement} are the results produced by the FC method from \cite{Lyon2010a}.  The errors for the FC method are reported using $\Delta x$ in \cite{Lyon2010}; in \Cref{fig:2D_refinement_heat_hard} we plot the results as a function of $n$ so that the values of $\Delta x$ are equivalent (accounting for the rescaling of the problem).  For both the IBSE and FC methods, the timestep has been taken low enough so that the dominant error is the spatial error.  Despite the fact that the FC method is fifth-order in space, the IBSE-$3$ method is able to achieve comparable errors at the finest grid-spacing reported\footnote{\color{review_color}The method used in \cite{Lyon2010} was based on an oversampled formulation near the boundary.  Newer FC methods \cite{Bruno2014} do not require this oversampling, effectively halving the number of discretization points in each dimension required to obtain a given accuracy.}.  For all of the FC results, the timestep is taken to be $1\times10^{-5}$.  Although their Alternating-Direction Implicit (ADI) scheme is unconditionally stable, it is only second-order accurate, and thus requires the use of small timesteps to obtain errors small enough to not impede the rapid convergence of the spatial discretization.  The ADI nature of their scheme complicates implementation of higher-order schemes, although this issue may be resolved using high-order Richardson extrapolation \cite{Albin2011}.  The structure of the IBSE method allows for the use of simple and high-order timestepping methods.  In these computations, we use BDF4 and are able to take substantially larger timesteps while still obtaining comparable error.

\section{Results: nonlinear problems}
\label{section:nonlinear}

One significant advantage of the IBSE method is the simplicity with which the elliptic and parabolic equation solvers can be integrated into methods to solve more complicated PDE.  We solve two nonlinear problems, the Burgers' equation and the Fitzhugh-Nagumo equations, as a demonstration.  The IBSE method is able to obtain high-order convergence in both space and time using simple timestepping methods and a straightforward pseudo-spectral computation of nonlinearities.

\subsection{Burgers' equation}
\label{section:nonlinear:burgers}

\begin{figure}
	\centering
	\hspace*{\fill}
	\begin{subfigure}[b]{0.4\textwidth}
		\centering
		\includegraphics[width=\textwidth]{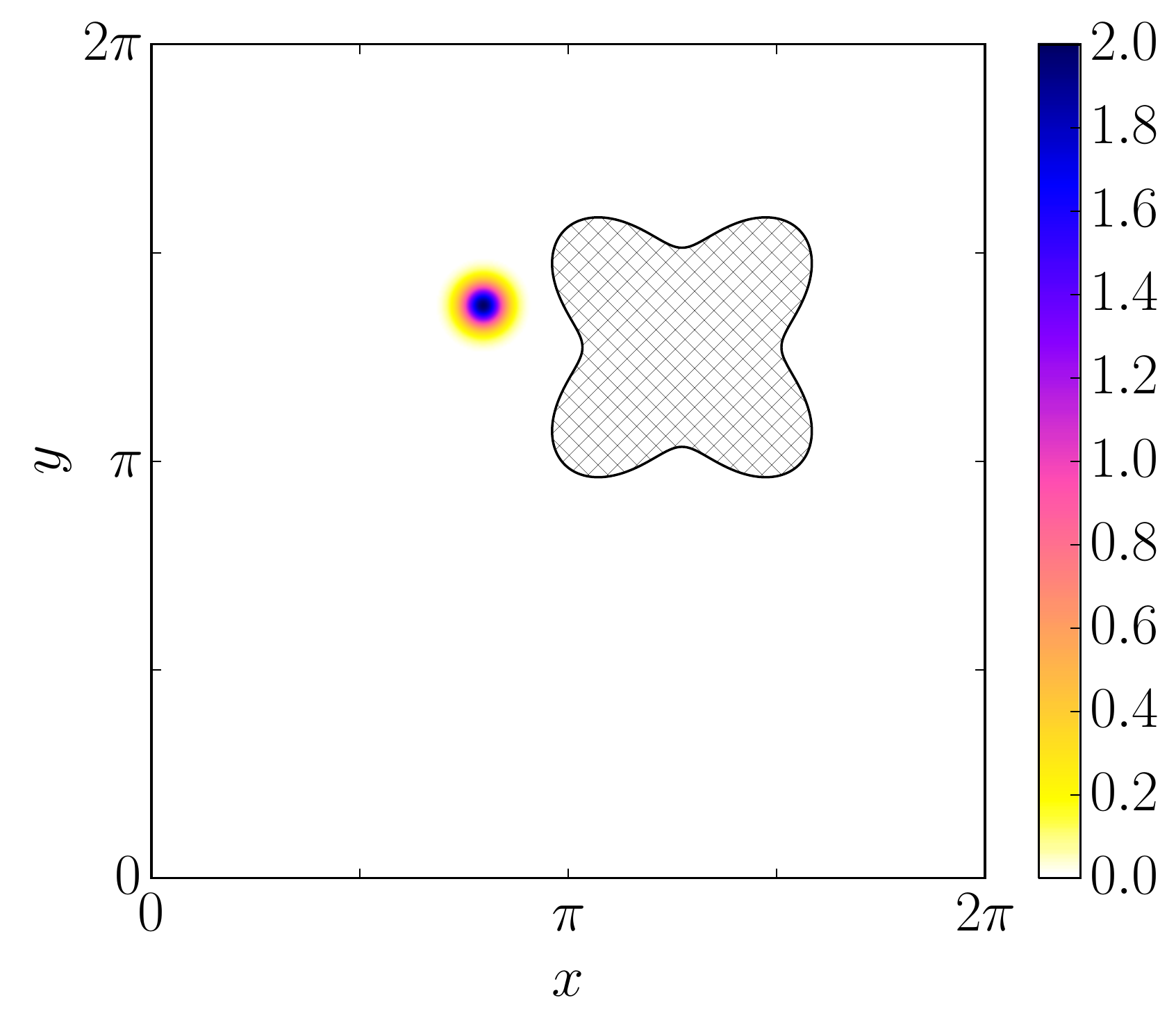}
		\subcaption{initial condition and geometry}
		\label{burgers:initial}
	\end{subfigure}
	\hfill
	\begin{subfigure}[b]{0.4\textwidth}
	\centering
		\includegraphics[width=\textwidth]{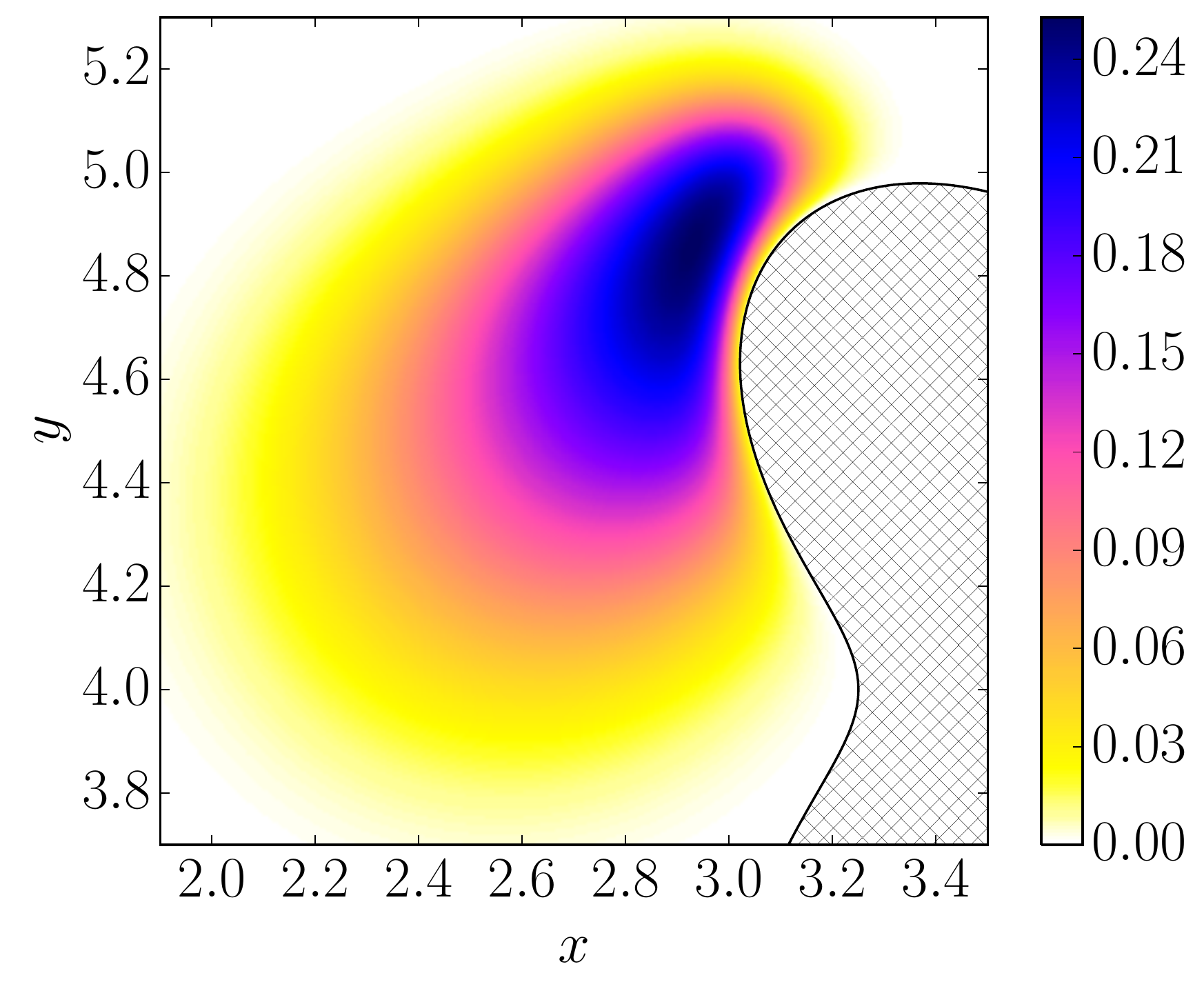}
		\subcaption{zoom of solution at $t=2$}
		\label{fig:burgers:solution}
	\end{subfigure}
	\hspace*{\fill}
	\caption{The initial condition, geometry, and solution at $t=2.0$ to \Cref{eq:burgers}.  The extension domain $E$ is indicated by the hatched region.  The physical domain $\Omega$ is everything outside of $E$.}
	\label{fig:burgers}
\end{figure}

In this section, we apply the IBSE method to a two-dimensional homogeneous Burgers' equation, a simple model with a nonlinearity that has the same form as Navier-Stokes:
\begin{subequations}
	\label{eq:burgers}
	\begin{align}
		u_t - \nu\Delta u + u\cdot\nabla u	&=	0	&	&\text{in }\Omega,	\\
		u									&=	0	&	&\text{on }\Gamma,	\\
		u(x,y,t=0)							&=	\psi(x,y),
	\end{align}
\end{subequations}
Here $\Omega$ is the region \emph{outside} of the polar region defined by $r=1+\cos(\theta+\pi/4)/4$, shifted by $(3\pi/2,3\pi/2)$; this domain is shown in \Cref{burgers:initial}.  The initial condition is given by the Gaussian pulse $\psi(x,y)=2e^{-40(x-2.5)^2}e^{-40(y-4.3)^2}$, and $\nu=0.01$.  As with the other examples, the computational domain $C$ is taken to be $\mathbb{T}^2=[0,2\pi]\times[0,2\pi]$.

We integrate this equation in time using a fourth-order implicit-explicit (IMEX) Backward Differentiation formula \cite{Hundsdorfer2007}:
\begin{multline}
	\frac{25}{12}u^{n+1} - 4u^n + 3u^{n-1} - \frac{4}{3}u^{n-2} + \frac{1}{4}u^{n-3} = \\
		\Delta t \left[\mathcal{I}(u^{n+1}) + 4\mathcal{E}(u^n) - 6\mathcal{E}(u^{n-1}) + 4\mathcal{E}(u^{n-2}) - \mathcal{E}(u^{n-3})\right],
\end{multline}
where $\mathcal{I}$ evaluates the \emph{stiff} terms in $u_t$ and $\mathcal{E}$ evaluates the \emph{non-stiff} terms in $u_t$.  For this problem, $\mathcal{I}(u)=\nu\Delta u$ and $\mathcal{E}(u)={\color{review_color}-}u\cdot\nabla u$.  In practice, timestepping involves the explicit computation of a forcing function:
\begin{equation}
	\frac{25}{12}f = 4u^n - 3u^{n-1} + \frac{4}{3}u^{n-2} - \frac{1}{4}u^{n-3} + \Delta t\left[4\mathcal{E}(u^n)-6\mathcal{E}(u^{n-1})+4\mathcal{E}(u^{n-2})-\mathcal{E}(u^{n-3})\right],
\end{equation}
followed by solution of the equation
\begin{subequations}
	\label{eq:burgers_update}
	\begin{align}
		\left(\mathbb{I}-\frac{12}{25}\nu\Delta t\Delta\right) u^{n+1} &= f	& &\text{in }\Omega,	\\
		u^{n+1}	&=	0	&	&\text{on }\Gamma.
	\end{align}
\end{subequations}
Since $u^{n}$ is smooth in the entire computational domain $C$, the explicit terms such as $\mathcal{E}(u^{n})$ may be accurately computed using simple pseudo-spectral methods.  \Cref{eq:burgers_update} is solved using the method described in \Cref{section:numerics:inversion}.  For this test, the startup values $u^{-1}$, $u^{-2}$, and $u^{-3}$ are all taken to be $\psi$.  We integrate \Cref{eq:burgers} to $t=2$ with a timestep of $\Delta t=\Delta x/20$, for $n=2^6$ to $n=2^{10}$, using IBSE-$3$ when solving \Cref{eq:burgers_update}.  \Cref{burgers:initial} and \Cref{fig:burgers:solution} show the initial condition and a zoom of the solution at $t=2$ for $n=2^{10}$.  Convergence is assessed by comparing the ratio of the $L^\infty$ difference between solutions at different resolutions, computed over the spatial locations in the coarsest grid.  Results are shown in \Cref{table:burgers_refinement}, demonstrating fourth-order convergence in both space and time.

\begin{table}
	\centering
	\begin{tabular}{cccc}
		\toprule
		Solutions			&	$L^\infty$ difference	&	ratio	&	$\log_2$ ratio	\\
		\midrule
		$2^6$, $2^7$		&	3.65e-02				&	-		&	-				\\
		$2^7$, $2^8$		&	9.75e-03				&	3.7		&	1.9				\\
		$2^8$, $2^9$		&	6.87e-04				&	14.2	&	3.8				\\
		$2^9$, $2^{10}$		&	2.83e-05				&	24.3	&	4.6				\\
		\bottomrule
	\end{tabular}
	\caption{Unnormalized $L^\infty$ differences between solutions to \Cref{eq:burgers} at successive levels of grid refinement, computed using IBSE-$3$.  Ratios of the $L^\infty$ differences are computed, indicating fourth-order convergence in space and time.}
	\label{table:burgers_refinement}
\end{table}

\subsection{Fitzhugh-Nagumo equation}
\label{section:nonlinear:fhn}

\begin{figure}
	\centering
	\hspace*{\fill}
	\begin{subfigure}[b]{0.32\textwidth}
		\centering
		\includegraphics[width=\textwidth]{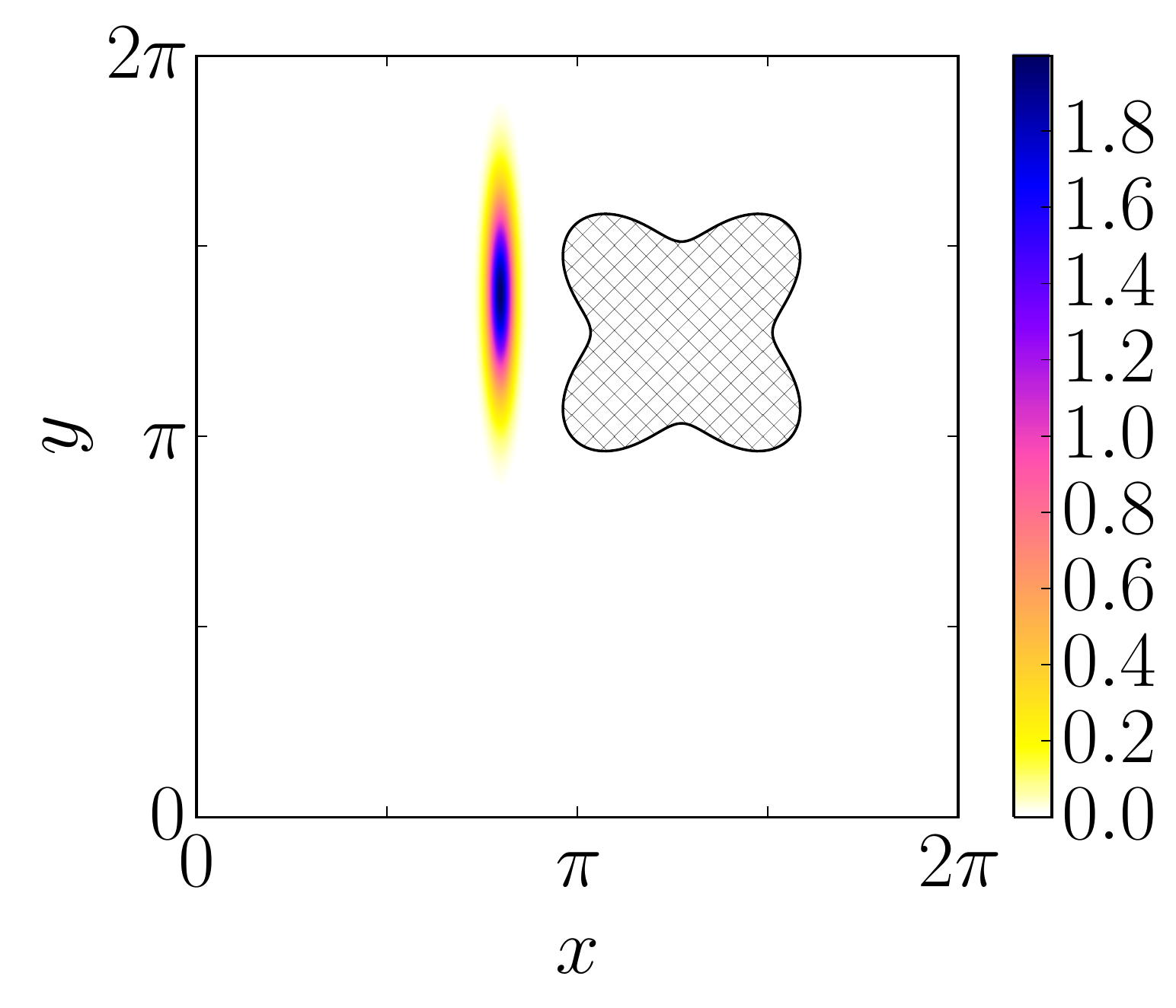}
		\subcaption{initial condition and geometry}
		\label{fig:fhn:initial}
	\end{subfigure}
	\hfill
	\begin{subfigure}[b]{0.32\textwidth}
	\centering
		\includegraphics[width=\textwidth]{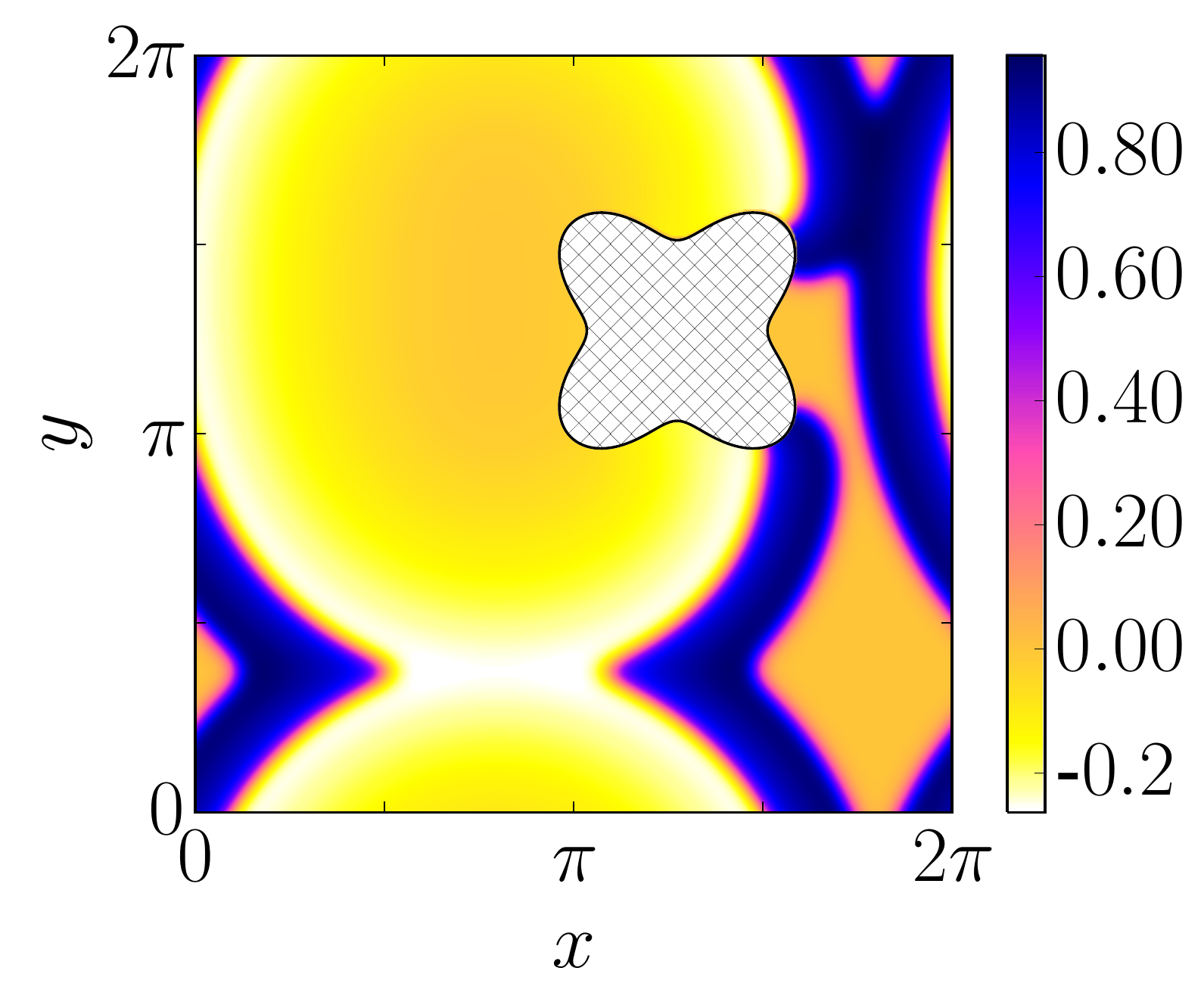}
		\subcaption{$v$ at $t=200$}
		\label{fig:fhn:solution:v}
	\end{subfigure}
	\hfill
	\begin{subfigure}[b]{0.32\textwidth}
	\centering
		\includegraphics[width=\textwidth]{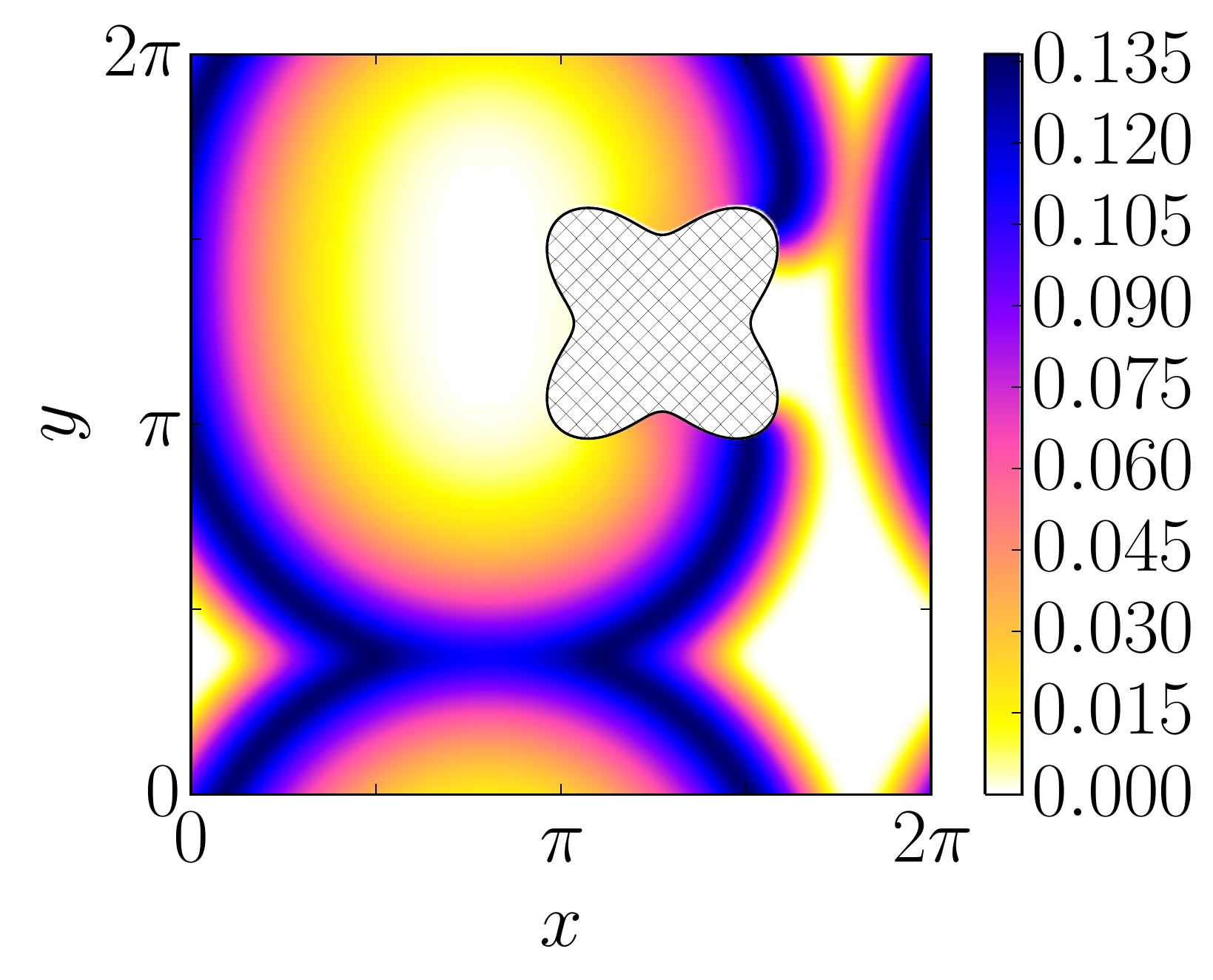}
		\subcaption{$w$ at $t=200$}
		\label{fig:fhn:solution:w}
	\end{subfigure}
	\hspace*{\fill}
	\caption{The initial condition for $v$, and the solutions to $v$ and $w$ at $t=200$ for \Cref{eq:fhn} generated with the IBSE-$3$ method.  The extension domain $E$ is indicated by the hatched region.  The physical domain $\Omega$ is everything outside of $E$.}
	\label{fig:fhn}
\end{figure}

The Fitzhugh-Nagumo equation \cite{Nagumo1962} is a system of nonlinear reaction-diffusion equations often used to model excitable media in biological systems (see e.g. \cite{aliev1996simple}).  We solve the Fitzhugh-Nagumo equations with homogeneous Neumann boundary conditions and an initially unperturbed recovery variable $w$:
\begin{subequations}
	\label{eq:fhn}
	\begin{align}
		v_t - \nu\Delta v + v(a-v)(1-v) + w		&=	0	&	&\text{in }\Omega,	\\
		w_t + \epsilon(\gamma w - v) 			&=	0	&	&\text{in }\Omega,	\\
		\frac{\partial v}{\partial n}			&=	0	&	&\text{on }\Gamma	\\
		v(x,y,t=0)	&=	\psi(x,y)												\\
		w(x,y,t=0)	&=	0.
	\end{align}
\end{subequations}
The physical domain $\Omega$ and the computational domain $C$ are the same as those used for the solution to Burgers' equation in \Cref{section:nonlinear:burgers}.  The initial data is the pulse $\psi(x,y)=2e^{-100(x-2.5)^2}e^{-10(y-4.3)^2}$.  High-order time integration is implemented using the same IMEX-BDF scheme used in \Cref{section:nonlinear:burgers}, with the implicit terms given by
\begin{equation}
	\mathcal{I}
	\begin{pmatrix}
		v	\\	w
	\end{pmatrix}
	=
	\begin{pmatrix}
		\nu\Delta v	\\
		0
	\end{pmatrix}
\end{equation}
and the explicit terms given by
\begin{equation}
	\mathcal{E}
	\begin{pmatrix}
		v	\\	w
	\end{pmatrix}
	=
	\begin{pmatrix}
		v(a-v)(v-1) - w	\\
		\epsilon(v-\gamma w).
	\end{pmatrix}
\end{equation}
We set the parameters to be $a=0.1$, $\gamma=2$, $\epsilon=0.005$, and $\nu=0.001$.  For this test, the startup values $v^{-1}$, $v^{-2}$, and $v^{-3}$ are all taken to be $\psi$.  We integrate \Cref{eq:fhn} to $t=200$ with a timestep of $\Delta t=\Delta x/2$, for $n=2^6$ to $n=2^{10}$ with the IBSE-$3$ method when solving for the updated value of $v$.  \Cref{fig:fhn} shows the initial condition $\psi$, as well as the solutions $v$ and $w$ at time $t=200$.  Convergence is assessed by comparing the ratio of the $L^\infty$ difference between solutions at different resolutions, computed over the spatial locations in the coarsest grid.  Results are shown in \Cref{table:fhn_refinement}, demonstrating at least third-order convergence in both space and time, as expected when solving with Neumann boundary conditions.

\begin{table}
	\centering
	\hspace*{\fill}
	\begin{subtable}[b]{0.45\textwidth}
		\begin{tabular}{cccc}
			\toprule
			Solutions			&	$L^\infty$ difference	&	ratio	&	$\log_2$ ratio	\\
			\midrule
			$2^6$, $2^7$		&	2.23e-01				&	-		&	-				\\
			$2^7$, $2^8$		&	2.38e-02				&	9.4		&	3.2				\\
			$2^8$, $2^9$		&	1.86e-03				&	12.8	&	3.7				\\
			$2^9$, $2^{10}$		&	2.68e-04				&	6.9		&	2.8				\\
			\bottomrule
		\end{tabular}
		\caption{Refinement study for $v$}
		\label{table:fhn_refinement:u}
	\end{subtable}
	\hfill
	\begin{subtable}[b]{0.45\textwidth}
		\begin{tabular}{cccc}
			\toprule
			Solutions			&	$L^\infty$ difference	&	ratio	&	$\log_2$ ratio	\\
			\midrule
			$2^6$, $2^7$		&	1.63e-02				&	-		&	-				\\
			$2^7$, $2^8$		&	2.25e-03				&	7.3		&	2.9				\\
			$2^8$, $2^9$		&	1.91e-04				&	11.8	&	3.6				\\
			$2^9$, $2^{10}$		&	8.68e-06				&	22.0	&	4.5				\\
			\bottomrule
		\end{tabular}
		\caption{Refinement study for $w$}
		\label{table:fhn_refinement:w}
	\end{subtable}
	\hspace*{\fill}
	\caption{Unnormalized $L^\infty$ differences between solutions to \Cref{eq:fhn} at successive levels of grid refinement, computed with IBSE-$3$.  Ratios of the $L^\infty$ differences are computed, indicating at least third-order convergence in space and time.}
	\label{table:fhn_refinement}
\end{table}



\section{Discussion}

We have developed a new numerical method, the IBSE method, for solving elliptic and parabolic PDE on general smooth domains to an arbitrarily high order of accuracy using simple Fourier spectral methods.  Accuracy is obtained by forcing the solution to be globally smooth on the entire computational domain.  This is accomplished by solving for the smooth extension of the unknown solution, enabling the method to directly solve elliptic equations.  Remarkably, this computation can be effectively reduced to a small, dense system of equations, and the cost to invert this system can be minimized by precomputing and storing its LU-factorization.  This provides a very efficient algorithm for implicit timestepping on stationary domains.  The method requires minimal geometric information regarding the boundary: only its position, normals, and an indicator variable denoting whether points in the computational grid are inside of the physical domain $\Omega$ or not, enabling simple and robust code.  The method is flexible enough to allow high-order discretization in space and time for a wide range of nonlinear PDE using straightforward implicit-explicit timestepping schemes.

The IBSE method shares some similarities with Fourier Continuation (FC) \cite{Lyon2010a,Lyon2010,Albin2011} and Active Penalty (AP) \cite{Shirokoff2013} methods.  The FC method uses Fourier methods to obtain a high-order discretization, and the idea that smooth, non-periodic functions can be turned into smooth and periodic functions by extending them into a larger domain in some appropriate way.  There are two key differences between the FC method and the IBSE method.  The FC method ({\small\emph{i}}) uses dimensional splitting to reduce the problem to a set of one-dimensional problems and ({\small\emph{ii}}) relies on data from the previous timestep in order to generate the Fourier continuations that allow for their high-order spatial accuracy.  This forces the FC method to use an Alternating-Direction Implicit scheme to take large stable timesteps when solving parabolic equations, complicating the implementation of high-order timestepping.  In addition, the FC method requires the use of an iterative solver for computing solutions to the Poisson equation, complicating an efficient discretization of the \emph{incompressible} Navier-Stokes equations, although the FC method has been used to solve the compressible Navier-Stokes equations to high order \cite{Albin2011}.  In contrast to the conditioning issues faced by the IBSE method, the one-dimensional nature of the Fourier-continuation problem allows the FC method to use high-precision arithmetic in certain precomputations to effectively eliminate the conditioning problem inherent in constructing smooth extensions \cite{Platte2011}.  This enables the FC method to obtain stable methods that converge to higher order than can be achieved by the IBSE method in double-precision arithmetic.

The AP method, like the FC method, relies on data from previous timesteps in the way that it imposes smoothness on the solution of the PDE.  A large drag force is applied that penalizes deviations of the solution in the extension region from the smooth extension of the solution at the previous timestep.  This dependence on data from previous timesteps forces the AP method to use \emph{explicit} timestepping when solving parabolic equations.  The smoothness constraint is also enforced only \emph{approximately}; solutions are $C^1$ regardless of how many derivatives are matched when generating the extension functions.  We believe this property explains both the relatively small difference in $L^\infty$ error between the second-order AP method and the second-order IBSE-$1$ method, and the orders of magnitude difference in $L^\infty$ error between the third-order AP method and the third-order IBSE-$2$ method observed in solutions to \Cref{eq:2D_test_heat}.  Despite these disadvantages, the AP method does not require the solution of a dense system of equations and should be immediately applicable to moving boundary problems.

We have only implemented two-dimensional examples in this paper.  The IBSE method extends to three dimensions without any changes, although there is one minor difficulty that must be resolved: the production of an accurate enough quadrature for the discretization of the boundary integrals in the spread operators $S$ and $T_k$.  For two-dimensional problems, this comes nearly for free since the simple quadrature rule given in \Cref{section:numerics:operators} is spectrally accurate for closed one-dimensional curves.  High-order quadrature rules for two-dimensional surfaces are more complicated but well-developed, and high-order surface representation has been incorporated into the Immersed Boundary method \cite{griffith2012hybrid,Shankar2012}.

We have discretized the IBSE-$k$ method for $1\leq k\leq 3$.  In this work, we are limited to $k=3$ largely by the \emph{accuracy and regularity of our regularization of $\delta$}.  This limitation is not fundamental.  There are at least two paths forward to obtain higher-order accuracy.  Smoother and more accurate $\delta$-functions with compact support could be constructed, however the method that we use to construct $\tilde\delta$ (convolving $\delta$-functions with limited regularity against themselves) would most likely not be of use here.  For example, to discretize the IBSE-$5$ method, we could start with the sixth-order analog to the base $\delta$-function that we use (see \Cref{section:numerics:operators:delta}), which has a support of six grid points \cite{Bringley2008}.  To obtain a $C^5$ $\delta$-function, this would need to be convolved against itself 5 times, resulting in a $\delta$ with a support of 30 gridpoints (900 for a two-dimensional $\delta$-function).  An alternate approach is to use non-compactly supported $\delta$-functions, combined with fast transform methods for the application of the spread and interpolation operators \cite{Greengard1991,Greengard2006}.

When solving the IBSE-$k$ equations for high values of $k$, the conditioning of the extension operator $\mathcal{H}^k$ given in \Cref{eq:extension_operator_definition} becomes problematic.  Although the conditioning can be controlled by an appropriate choice of $\Theta$, the length scales introduced by $\Theta$ rapidly become unresolvable as the discretization is refined.  A different choice of extension operator, e.g. a polynomial of the Laplacian whose spectra is explicitly designed to minimize the condition number while limiting the introduction of fine scales, could potentially be used to mitigate this problem.

In this manuscript, we have developed two of the necessary components for the implementation of a high-order numerical scheme for the incompressible Navier-Stokes equations based on the IBSE method: a direct solver for the Poisson equation and a high-order implicit-explicit discretization of a nonlinear advection-diffusion equation.  The missing component is \emph{how to impose a global divergence constraint} without compromising the accuracy of the discretization.  This is a non-trivial problem even for finite-difference discretizations on curvilinear domains \cite{Shirokoff2010}.  We have obtained preliminary results for solving the Stokes equation using the IBSE method; this work will be presented in a forthcoming contribution.

Finally, the high-efficiency that we obtain is only currently achievable on stationary domains.  This is because the inversion method that we use in this paper for the IBSE method relies on an expensive precomputation that depends on the physical domain and the discretization.  However, there is no fundamental obstacle to the application of the IBSE method to problems with moving domains.  Instead, the challenge is to find a robust and efficient method to invert the IBSE system in \Cref{eq:the_system} that does not require substantial precomputation.  Recent progress has been made for preconditioning similar systems of equations for the simulation of rigid-body motion in an Immersed Boundary framework \cite{kallemov2015immersed}.  The integration of these ideas into the IBSE method to allow for simulation of moving boundary problems will be an area of active future research.



\section*{Acknowledgements}

This work was supported in part by the National Science Foundation under Grant DMS-1160438.  The authors would like to thank Grady Wright for useful discussions regarding numerical conditioning problems that helped lead to the form of the operator $\mathcal{H}^k$ and the definition of $\Theta^*$ in \Cref{section:numerics:extension}.



\clearpage

\bibliography{bib,manual_bib}{}
\bibliographystyle{unsrtnat}


\clearpage

\begin{appendices}
\crefalias{section}{appsec}
\crefalias{subsection}{appsec}

\section{Formula for $C^3$ $\delta$-function with interpolation accuracy of $\mathcal{O}(\Delta x^4)$}
\label{appendix:delta}

We give the formula here for $\tilde\delta$, defined in \Cref{section:numerics:operators:delta}.  For brevity, we define $\tilde\delta$ for $r\geq 0$, which fully defines the function as $\tilde\delta$ is even.  The formula for $\tilde\delta$ is defined piecewise; each portion is a fifteenth order polynomial.  Coefficients are provided in the table below.

\vspace{1em}

{
	\centering
	\setlength{\tabcolsep}{2pt}

	\renewcommand{\arraystretch}{1.75}

	\begin{tabular}{c|c|c|c|c|c|c|c|c|c|}
		&	$0 \leq r < 1$	&	$1 \leq r < 2$	&	$2 \leq r < 3$	&	$3 \leq r < 4$	&	$4 \leq r < 5$	&	$5 \leq r < 6$	&	$6 \leq r < 7$	&	$7 \leq r < 8$	&	$r >= 8$	\\
		\hline\hline
		$r^{0}$	&	$\frac{12949745023}{20432412000}$	&	$\frac{3177441629}{5003856000}$	&	$\frac{21914742667}{35026992000}$	&	$\frac{4094824493}{17513496000}$	&	$\frac{-1606651889}{5837832000}$	&	$\frac{163201885541}{35026992000}$	&	$\frac{1005507698627}{245188944000}$	&	$\frac{-5005877248}{1915538625}$	&0	\\
		\hline
		$r^{1}$	&	$0$	&	$\frac{-171811}{22861440}$	&	$\frac{-2566373}{38102400}$	&	$\frac{3468455}{4191264}$	&	$\frac{301286857}{62868960}$	&	$\frac{-4590637187}{1089728640}$	&	$\frac{-14398288259}{1089728640}$	&	$\frac{421762048}{127702575}$	&0	\\
		\hline
		$r^{2}$	&	$\frac{-16459}{30240}$	&	$\frac{-9089}{17280}$	&	$\frac{-14339}{51840}$	&	$\frac{-163307}{285120}$	&	$\frac{-703993}{95040}$	&	$\frac{-2987689}{673920}$	&	$\frac{56979607}{3991680}$	&	$\frac{-23168}{45045}$	&0	\\
		\hline
		$r^{3}$	&	$0$	&	$\frac{-64649}{3265920}$	&	$\frac{-1946191}{5443200}$	&	$\frac{-3764137}{2993760}$	&	$\frac{32748677}{8981280}$	&	$\frac{101232611}{11975040}$	&	$\frac{-17097977}{2395008}$	&	$\frac{-2091088}{1403325}$	&0	\\
		\hline
		$r^{4}$	&	$\frac{81491}{453600}$	&	$\frac{141751}{777600}$	&	$\frac{288599}{777600}$	&	$\frac{6701891}{4276800}$	&	$\frac{560257}{1425600}$	&	$\frac{-4187303}{777600}$	&	$\frac{78901349}{59875200}$	&	$\frac{642734}{467775}$	&0	\\
		\hline
		$r^{5}$	&	$0$	&	$\frac{88517}{5443200}$	&	$\frac{61633}{3024000}$	&	$\frac{-93301}{151200}$	&	$\frac{-1620853}{1360800}$	&	$\frac{1038857}{604800}$	&	$\frac{255833}{604800}$	&	$\frac{-538927}{850500}$	&0	\\
		\hline
		$r^{6}$	&	$\frac{-11737}{340200}$	&	$\frac{-119603}{2332800}$	&	$\frac{-269467}{2332800}$	&	$\frac{16957}{1166400}$	&	$\frac{218939}{388800}$	&	$\frac{-488741}{2332800}$	&	$\frac{-5818427}{16329600}$	&	$\frac{773411}{4082400}$	&0	\\
		\hline
		$r^{7}$	&	$\frac{143}{145152}$	&	$\frac{298727}{45722880}$	&	$\frac{630773}{15240960}$	&	$\frac{1057927}{15240960}$	&	$\frac{-1137851}{9144576}$	&	$\frac{-152395}{3048192}$	&	$\frac{1823393}{15240960}$	&	$\frac{-1825543}{45722880}$	&0	\\
		\hline
		$r^{8}$	&	$\frac{3223}{793800}$	&	$\frac{28127}{5443200}$	&	$\frac{-7337}{5443200}$	&	$\frac{-9859}{388800}$	&	$\frac{7069}{907200}$	&	$\frac{157289}{5443200}$	&	$\frac{-966457}{38102400}$	&	$\frac{58621}{9525600}$	&0	\\
		\hline
		$r^{9}$	&	$\frac{-143}{435456}$	&	$\frac{-30173}{19595520}$	&	$\frac{-79651}{32659200}$	&	$\frac{4771}{1306368}$	&	$\frac{61009}{19595520}$	&	$\frac{-44227}{6531840}$	&	$\frac{24421}{6531840}$	&	$\frac{-69019}{97977600}$	&0	\\
		\hline
		$r^{10}$	&	$\frac{-3211}{13608000}$	&	$\frac{-403}{11664000}$	&	$\frac{7033}{11664000}$	&	$\frac{91}{2916000}$	&	$\frac{-247}{243000}$	&	$\frac{11519}{11664000}$	&	$\frac{-32227}{81648000}$	&	$\frac{2447}{40824000}$	&0	\\
		\hline
		$r^{11}$	&	$\frac{13}{483840}$	&	$\frac{13}{207360}$	&	$\frac{-13}{345600}$	&	$\frac{-221}{2280960}$	&	$\frac{13}{84480}$	&	$\frac{-221}{2280960}$	&	$\frac{53}{1774080}$	&	$\frac{-299}{79833600}$	&0	\\
		\hline
		$r^{12}$	&	$\frac{1}{181440}$	&	$\frac{-1}{155520}$	&	$\frac{-1}{155520}$	&	$\frac{7}{427680}$	&	$\frac{-1}{71280}$	&	$\frac{1}{155520}$	&	$\frac{-19}{11975040}$	&	$\frac{1}{5987520}$	&0	\\
		\hline
		$r^{13}$	&	$\frac{-1}{1451520}$	&	$\frac{-1}{2612736}$	&	$\frac{29}{21772800}$	&	$\frac{-13}{9580032}$	&	$\frac{113}{143700480}$	&	$\frac{-173}{622702080}$	&	$\frac{1}{17791488}$	&	$\frac{-47}{9340531200}$	&0	\\
		\hline
		$r^{14}$	&	$\frac{-1}{25401600}$	&	$\frac{1}{10886400}$	&	$\frac{-1}{10886400}$	&	$\frac{1}{17107200}$	&	$\frac{-1}{39916800}$	&	$\frac{1}{141523200}$	&	$\frac{-1}{838252800}$	&	$\frac{1}{10897286400}$	&0	\\
		\hline
		$r^{15}$	&	$\frac{1}{203212800}$	&	$\frac{-1}{261273600}$	&	$\frac{1}{435456000}$	&	$\frac{-1}{958003200}$	&	$\frac{1}{2874009600}$	&	$\frac{-1}{12454041600}$	&	$\frac{1}{87178291200}$	&	$\frac{-1}{1307674368000}$	&0	\\
		\hline
		\hline
	\end{tabular}

	\renewcommand{\arraystretch}{1}
}

\section{Inversion of the IBSE-$k$ system \eqref{eq:the_system} for $\mathcal{L}=\Delta$}
\label{appendix:poisson_inversion}

For the special case of a Fourier discretization of the Poisson problem, where $\mathcal{L}$ is the periodic Laplacian, the method of inversion given in \Cref{section:numerics:inversion} is complicated by the nullspace of $\Delta$.  In particular, the main block in \Cref{eq:the_system},
\begin{equation}
	\begin{pmatrix}
		\Delta	&	-\chi_E\Delta	\\
				&	\mathcal{H}^k
	\end{pmatrix}
\end{equation}
is not invertible, and thus the Schur-complement $SC$ \eqref{eq:the_schur_complement} cannot be directly formed and prefactored.  The resolution of this problem is conceptually simple but computationally involved, so we will give a detailed description for the simpler case of a direct-forcing Immersed Boundary discretization and state the result for our problem.  A direct-forcing IB discretization of the Poisson problem with Dirichlet boundary conditions is
\begin{equation}
	\begin{pmatrix}
		\Delta	&	S	\\
		S^*
	\end{pmatrix}
	\begin{pmatrix}
		u	\\	G
	\end{pmatrix}
	\begin{pmatrix}
		f	\\	b
	\end{pmatrix}.
	\label{eq:discrete_IB_poisson}
\end{equation}
Formally, $G$ may be computed by solving
\begin{equation}
	\label{eq:IB_schur}
	(S^*\Delta^{-1}S)G = S^*\Delta^{-1}f - b.
\end{equation}
Despite the fact that \Cref{eq:discrete_IB_poisson} is invertible, the periodic Laplacian has a nullspace, making direct implementation of this approach to compute $G$ impossible.  Instead, let us decompose the solution $u$ as
\begin{equation}
	u = u_0 - Uv, \label{eq:decomposition}
\end{equation}
where $u_0$ has mean $0$, $U$ is a scalar, and $v=\mathds{1}(\Delta x)^d$, where $\mathds{1}$ denotes the vector of all ones.  We note that $v$ spans the nullspace of the self-adjoint operator $\Delta$, and we have scaled $v$ so that $v^\intercal u$ is equal to the the discrete integral of $u$.  We will also denote the averaging operator $\fint_C w=\abs{C}^{-1}\int_C w$, discretely this is given for $C$ by $\fint_C w = (2\pi)^{-d}v^\intercal w$ and for $\Gamma$ by $\fint_\Gamma H=(2\pi)^{-1}\Delta s^\intercal H$, where $\Delta s$ is the vector of quadrature weights on the boundary (see \Cref{section:numerics:operators}).  Plugging the decomposition given in \Cref{eq:decomposition} into \Cref{eq:discrete_IB_poisson}, we find that
\begin{equation}
	\label{eq:projected_poisson}
	\Delta u_0 + SG = f.
\end{equation}
Define a family of pseudo-inverses to $\Delta$, $\{\mathcal{A}_\mu\}_{\mu> 0}$, by its Fourier series to be
\begin{equation}
	\widehat{\mathcal{A}_\mu f} =
	\begin{cases}
		\mu\widehat{f_0}					&\abs{k} = 0,	\\
		\frac{-1}{\abs{k}^2}\widehat{f_k}	&\abs{k}\neq0.
	\end{cases}
\end{equation}
Because $u_0$ has mean $0$, $\mathcal{A}_\mu\Delta u_0=u_0$ for any $\mu$. Applying $\mathcal{A_\mu}$ to \Cref{eq:projected_poisson}, we find that $u_0 = \mathcal{A}_\mu f - \mathcal{A}_\mu SG$, and applying $S^*$ to this equation, along with the constraint that $S^* u=b$, yields that
\begin{equation}
	S^*\mathcal{A}_\mu SG + S^*(Uv) = S^*\mathcal{A}_\mu f - b,
\end{equation}
which is a reformulation of the formal equation \eqref{eq:IB_schur}.  This is now an equation for the \emph{two} unknowns $G$ and $U$, and must be supplemented with the additional equation found by projecting $\Delta u + SG = f$ onto the nullspace of $\Delta$, giving $v^\intercal SG = v^\intercal f$.  The combination of these two equations forms an augmented Schur-complement system
\begin{equation}
	\begin{pmatrix}
		S^*\mathcal{A}_\mu S	&	S^*v \\
		v^\intercal S			&
	\end{pmatrix}
	\begin{pmatrix}
		G	\\	U
	\end{pmatrix}
	=
	\begin{pmatrix}
		S^*\mathcal{A}_\mu f - b	\\	v^\intercal f
	\end{pmatrix}
\end{equation}
that may be directly formed and prefactored.  With $(G,U)$ known, we may compute $u_0=\mathcal{A}_\mu f - \mathcal{A}_\mu SG$, and finally $u=u_0-U v$.  This system may be simplified by exploiting the fact that the nullspace of $\Delta$ is only constant functions and that $\fint_C SG=\fint_\Gamma G$, giving the equivalent system
\begin{equation}
	\label{eq:augmented_SC_IB}
	\begin{pmatrix}
		S^*\mathcal{A}_\mu S	&	\beta^{-1}\mathds{1}	\\
		\Delta s^\intercal/2\pi
	\end{pmatrix}
	\begin{pmatrix}
		G	\\	U
	\end{pmatrix}
	=
	\begin{pmatrix}
		S^*\mathcal{A}_\mu f - b	\\
		\fint_C f
	\end{pmatrix},
\end{equation}
where $\beta$ is a free parameter.  We remark that when $\Delta s$ is a constant vector, choosing $\beta=n_\text{bdy}$ symmetrizes the matrix in \Cref{eq:augmented_SC_IB}.  Although any choice of $\mu$ and $\beta$ may be used, we have observed empirically that $\mu=0$ and $\beta=n_\text{bdy}$ yields the most well-conditioned system.  Once $G$ and $U$ are known, we may compute $u$ to be
\begin{equation}
	u = \mathcal{A}_\mu(f - SG) - \beta^{-1} U.
\end{equation}
Following the same procedure for \Cref{eq:the_system} gives the augmented Schur-complement:
\begin{equation}
	\label{eq:augmented_SC}
	\begin{pmatrix}
		T_k^*(\mathcal{A}_\mu\chi_E\Delta-I)(\mathcal{H}^k)^{-1}T_k	&	T_k^*\mathcal{A}_\mu S		&	\beta^{-1}\mathbb{Y}_k	\\
		S^*\mathcal{A}_\mu\chi_E\Delta(\mathcal{H}^k)^{-1}T_k			&	S^*\mathcal{A}_\mu S	&	\beta^{-1}\mathds{1}	\\
		\fint_C\chi_E\Delta (\mathcal{H}^k)^{-1}T_k	&	\fint_C S
	\end{pmatrix}
	\begin{pmatrix}
		F	\\	G	\\	\lambda
	\end{pmatrix}
	=
	\begin{pmatrix}
		S^*\mathcal{A}_\mu\chi_\Omega f - g \\ T_k^*\mathcal{A}_\mu\chi_\Omega f \\ \fint_C\chi_\Omega f
	\end{pmatrix}.
\end{equation}
Here $\mathbb{Y}_k$ is defined by $
	\mathbb{Y}_k = 
	\begin{pmatrix}
		\mathds{1}_{n_\text{bdy}}	&	\mathds{O}_{(k-1)n_\text{bdy}}
	\end{pmatrix}^\intercal$;
this form is due to the fact that the estimates of \emph{values} produced by $T_k^*$ are affected by changes in mean while estimates of \emph{normal derivatives} are not.  Again, we empirically observe that choosing $\mu=0$ and $\beta=n_\text{bdy}$ provides the most well-conditioned system.  Once $F$ and $G$ are known, we can compute $\xi$ as
\begin{equation}
	\xi = -(\mathcal{H}^k)^{-1}T_kF,
\end{equation}
and once $\xi$ is known, $u$ can be computed as
\begin{equation}
	\label{eq:get_u_noninvertible}
	u = \mathcal{A}_\mu(\chi_D f + \chi_E\Delta\xi - SG) - \beta^{-1}U.
\end{equation}

\end{appendices}

\end{document}